\documentclass{amsart}
\usepackage{amsfonts} 
\usepackage{amsmath}
\usepackage{amssymb}
\usepackage{enumerate}
\usepackage{setspace} 
\usepackage{graphicx}
\newtheorem{theorem}{Theorem}[section]
\newtheorem{corollary}[theorem]{Corollary}
\newtheorem{lemma}[theorem]{Lemma}
\newtheorem{proposition}[theorem]{Proposition}
\theoremstyle{definition}
\newtheorem{definition}[theorem]{Definition}
\newtheorem{example}[theorem]{Example}

\newtheorem{remark}[theorem]{Remark}

\theoremstyle{remark}
\newtheorem*{acknowledgements}{Acknowledgements}
\numberwithin{equation}{section}

\begin{document}
\title[Asymptotics at infinity for fractal measures]{Affine systems:  Asymptotics at infinity for fractal measures}
\author[P.E.T. Jorgensen]{Palle E. T. Jorgensen}
\address[Palle E.T. Jorgensen]{Department of Mathematics, The University of Iowa, Iowa
City, IA 52242-1419, U.S.A.}
\email{jorgen@math.uiowa.edu}
\urladdr{http://www.math.uiowa.edu/\symbol{126}jorgen/}
\author[K.A. Kornelson]{Keri A. Kornelson}
\author[K.L. Shuman]{Karen L. Shuman}
\address[Keri Kornelson, Karen Shuman]{Department of Mathematics and Statistics,
Grinnell College, Grinnell, IA 50112-1690, U.S.A.}
\urladdr{http://www.math.grinnell.edu/\symbol{126}kornelso/}
\urladdr{http://www.math.grinnell.edu/\symbol{126}shumank/}
%
%
%
%
%
%
\thanks{This material is based upon work partially supported by the
U.S. National Science Foundation under grants DMS-0457581 and DMS-0503990, by the University of Iowa Department of Mathematics
NSF VIGRE grant DMS-0602242, and by the Grinnell College Committee for Support of Faculty Scholarship}

\begin{abstract}
     We study measures on $\mathbb{R}^d$ which are induced by a class of infinite
and recursive iterations in symbolic dynamics.  Beginning with a finite
set of data, we analyze prescribed recursive iteration systems, each
involving subdivisions. The construction includes measures arising from
affine and contractive iterated function systems with and without overlap (IFSs), i.e., limit
measures $\mu$  induced by a finite family of affine mappings in $\mathbb{R}^d$ (the
focus of our paper), as well as equilibrium measures in complex
dynamics.

     By a systematic analysis of the Fourier transform of the measure
$\mu$ at hand (frequency domain), we identify asymptotic laws, spectral
types, dichotomy, and chaos laws. In particular we show that the cases
when $\mu$ is singular carry a gradation, ranging from Cantor-like
fractal measures to measures exhibiting chaos, i.e., a situation when
small changes in the initial data produce large fluctuations in the
outcome, or rather, the iteration limit (in this case the measures). Our
method depends on asymptotic estimates on the Fourier transform of $\mu$
for paths at infinity in $\mathbb{R}^d$. We show how properties of $\mu$ depend on
perturbations of the initial data, e.g., variations in a prescribed
finite set of affine mappings in $\mathbb{R}^d$, in parameters of a rational
function in one complex variable (Julia sets and equilibrium measures),
or in the entries of a given infinite positive definite matrix.
\end{abstract}
\maketitle
\section{Introduction}
 We study measures associated with certain dynamical systems as
infinite iteration. Beginning with a finite set of data, we analyze
prescribed recursive iteration systems, each involving subdivisions.
Passing to the limit, we arrive at certain measures $\mu$ which capture
intriguing dynamical information, such as spectral type, and the
presence of ``chaos.''  By chaos, we refer to a situation when small
changes in the initial data produce large fluctuations in the outcome,
or rather, the iteration limit (in this case the measures). The initial
data we study may be a prescribed finite set of affine mappings in $\mathbb{R}^d$
(affine iterated function systems (IFSs)), it may be a rational
function in one complex variable (Julia sets and equilibrium measures),
or simply an infinite positive definite matrix (determinantal
measures).

       By a systematic analysis of the Fourier transform of the
measure $\mu$ at hand (frequency domain), we identify asymptotic laws,
spectral types, dichotomy, and chaos laws. Since our approach is
somewhat interdisciplinary we include details from one area that may
perhaps not be well know to readers from another. We are motivated by
several related developments and in particular by pioneering ideas of
Erd\H{o}s, Jessen, Kahane, and Wintner.

         Consider a fixed dimension $d$, and consider a given finite
set $S$ of affine and contractive mappings $\tau_i : \mathbb{R}^d \rightarrow \mathbb{R}^d$.  Iterations of the maps in $S$ generate what is called an ``affine Iterated Function System," abbreviated IFS.    There are several interesting limits in the literature associated with $S$; in this paper we will concentrate on two of these limits.  The first is an attractor $X(S)$, a compact subset of $\mathbb{R}^d$ which arises from recursive iterations of the maps in $S$.  Second, given a set of positive probabilities $p = \{p_i\}$ associated with the maps $\tau_i$ in $S$, there is a unique measure supported on $X(S)$ called the Hutchinson or equilibrium measure $\mu_{S,p}$.  We will study both the measure $\mu_{S,p}$ and generalizations of $\mu_{S,p}$.  The generalizations are induced by determinantal measures on infinite product spaces. 

Instances where a given system $(X(S), \mu_{S,p} )$ is known to be
singular relative to $d$-dimensional Lebesgue measure include the
Cantor-like fractals (e.g., $d = 1$, and the familiar deleted
middle-third construction). More generally, these are cases when $X(S)$
is ``thin'' as a subset of $\mathbb{R}^d$, when the associated fractal dimension
$d_S$ of $\mu_{S,p}$ is computable and satisfies  $d_S < d$. These systems
come with an intrinsic scale parameter (see e.g., (\ref{tauzeroone}); or a
scaling matrix (\ref{tauB})).
      For IFSs as the scale parameter expands, the gaps in the support
set $X(S)$ will ``close up,'' for example when the IFS-recursion has
overlaps. Then typically $X(S)$ has non-empty interior as a subset in
$\mathbb{R}^d$ and the properties of $\mu_{S,p}$  become more elusive. A very
special case of this was studied in 1939 for $d = 1$, and for a
single-variable scale parameter $\lambda$, in a paper by Erd\H{o}s \cite{Erd39}.
The occurence of these singular measures is somewhat
counter-intuitive, different from the Cantor fractals; we wish to
extract some of the underlying geometry of these singular measures.

Motivated by Erd\H{o}s's work on singular infinite Bernoulli convolution measures, we study which regularity and geometric properties might be gleaned from the Fourier transform of the Hutchinson measure $\mu_{S,p}$.   In particular, given $S$ and $p$, we are interested in locating measures $\mu_{S,p}$ which are \textit{not} relatively absolutely continuous with respect to the ambient Lebesgue measure.  In Section \ref{Sec:Pisot} we explain Erd\H{o}s's work for certain IFSs on $\mathbb{R}$, and in Sections \ref{Sec:2D} and \ref{Sec:higherdim} we turn to higher dimensions.

        In Section \ref{Sec:InducedMeasures}, we introduce  an important class of measures associated with IFSs which are
induced from
measures on infinite products. Since an IFS is a finite family $S$ of
mappings in an ambient space, for example in $\mathbb{R}^d$, we can
understand infinite iterations via the compact product space $P_S :=
S^\mathbb{N}$, or $S^\mathbb{Z}$, where $\mathbb{N} = \{1, 2,
\ldots\}$ denotes the natural numbers, and $\mathbb{Z}$ denotes the
integers.  

     In the past decade the literature on IFSs (self-affine sets and measures), encoding and digit representations for radix matrices has grown; in
part this growth is due to applications to such areas as number theory, dynamics,
and combinatorial geometry.   It is not possible here to give a complete
list of these directions.   Our present work has been motivated by the
following papers: \cite{AkSc05, Cur06, GaYu06, HLR02, HeLa04, KLSW99, Li06, Li07, Saf98, ZZ06, LaWa00}.  Here we are especially motivated by questions in number theory,
starting with Jessen-Wintner \cite{JeWi35} and Knuth \cite{Knu75}; see also \cite{AGPT06,
AlCe07}.

       Knuth's analysis \cite{Knu75} of what is now called affine IFSs was motivated by the desire to introduce geometry
into algorithms for general digit sets in positional number systems.  The
idea of turning ``digits'' into geometry and probability theory was followed
up later by others, e.g., Odlyzko \cite{Odl78}, Hutchinson \cite{Hut81},
Bratteli-Jorgensen \cite{BJ99}, and Lagarias-Wang \cite{LaWa96a, LaWa96b, LaWa96c, LaWa97, LaWa00}. 

      In fact Knuth \cite{Knu75} stresses that for a fixed radix number $a$, there
are many choices of digit sets $B$ which yield a positional number system. Similarly Knuth suggested the use of a matrix
in place of $a$. Rather than having a radix number and a finite set of
integers for digits in a positional number system, it was suggested that the
radix should be an expansive square matrix $A$ over the integers, i.e., some $A$
in $M_d(\mathbb{Z})$ for fixed $d$, and that the digit set $B$ should be a suitably chosen finite subset in the rank-$d$ lattice $\mathbb{Z}^d$.

       We shall adopt this setting, and our associated measures $\mu =
\mu_{(A,B)}$ (Section \ref{Subsec:IFSGeneralCase}) will then be supported on compact attractors $X = X(A,B)$ in $\mathbb{R}^d$, studied first in \cite{Hut81}.

        The aim of this paper is to revisit the Fourier asymptotics of the
measures $\mu_{(A,B)}$ in light of recent results on IFS involving dynamics and
representation theory, see e.g., \cite{DuJo06a, DuJo06b, DuJo06c, DuJo07b}.   The measures $\mu_{(A,B)}$ are equilibrium measures for the corresponding
affine system. In Section \ref{Sec:Pisot} we consider equilibrium measures for complex
iteration systems in the sense of \cite{Bea91}; in Section \ref{Sec:InducedMeasures}, we consider the analogous
measures in $\mathbb{R}^d$ induced by infinite determinants.

Section \ref{Sec:InducedMeasures} accomplishes four things.
First, we show how tools from operator theory serve to construct determinantal measures $\{\mu\}$ on $P_S$.  Then we show how for each pair
$(\mu, S)$, where $\mu$ is one of the measures and $S$ is a given IFS,
there is an induced measure on $\mathbb{R}^d$; and in Lemma
\ref{Lemma:BorelFourier} we give a formula for the Fourier transform
of the induced measures. Finally, we show how the larger class of
induced measures generalizes the class of Hutchinson measures
associated with IFSs.  We end with a discussion of how our Fourier asymptotics give rise to functions which we call infinite determinants.

\section{Background for affine IFSs in $\mathbb{R}^d$}\label{Sec:BackgroundIFS}
          
\subsection{One dimension:  $d=1$}

A special class of IFSs is a
one-parameter family of probability measures with compact support on
the real line $\mathbb{R}$ in which the probability measures
$\mu_{\lambda}$ arise as infinite convolutions of Bernoulli measures.
If $\lambda$ is fixed in the open unit interval $(0, 1)$, the
measure $\mu_\lambda$ is the Hutchinson measure \cite{Hut81} of the
simplest IFS on $\mathbb{R}$ given by  
\begin{equation}
\tau_0(x) = \lambda x \quad \text{ and } \quad \tau_1(x) = \lambda (x + 1).
\end{equation}\label{tauzeroone}%

\begin{example} A simple but illustrative case occurs when $\lambda = 1/3$---in this case the attractor $X(S)$ for $S$ in (\ref{tauzeroone}) is essentially Cantor's set with gaps. To construct Cantor's set, we make successive subdivisions of an
initially chosen interval using the two mappings $\tau_0$ and $\tau_1$
in (\ref{tauzeroone}) at each iteration step.  Also at each step, we
normalize the measure.  In the limit we arrive at Cantor's measure
$\mu_{1/3}$, the measure whose cumulative density is sometimes known as ``The
Devil's Staircase'' (see \cite[p. 38]{Fol84}).  This
measure $\mu_{1/3}$ is fractal in the strongest sense: its scaling
dimension $s$ is $s =\log_3(2)$ \cite{Hut81}.  Perhaps the more
familiar version of Cantor's construction is the one that leaves gaps
in ``the middle,'' i.e., the one associated with the two mappings $x
\mapsto x/3$ and $x \mapsto (x + 2)/3$ (see \cite[p. 74]{Jor06}).  However, the choice of where the gaps are placed
doesn't affect our present consideration; the presence of the gaps and
the relative size of the gaps are what is essential.
\end{example}\label{Ex:Cantor1/3Intro}%

In (\ref{tauzeroone}), as the value of $\lambda$ increases through the
interval $[1/3, 1/2)$, the measures $\mu_{\lambda}$ remain fractal and purely singular (the Hausdorff dimension is strictly less than $1$), but
the scaling dimension $s$ depends on $\lambda$; specifically, $s = -
\log(2)/\log(\lambda)$.  When $\lambda = 1/2$, the measure $\mu_{1/2}$
is Lebesgue measure restricted to an interval.  As $\lambda$ increases
beyond $1/2$, the subdivisions resulting from iteration of the system
(\ref{tauzeroone}) in 1D create overlap, so one might expect that the measures in the iteration limit
would become ``nicer'' and certainly not to have fractal features.  In addition, in the overlap case, no open set condition is satisfied, and the scaling dimension is not known. 

 Erd\H{o}s \cite{Erd39} proved that this reflection symmetry is broken for $d=1$: the situation
is more complicated.  In fact, the theorem of Erd\H{o}s \cite{Erd39}
states that the infinite Bernoulli convolution measure $\mu_{\lambda}$
is purely singular if $\lambda$ is the reciprocal of a Pisot number.
We will discuss Pisot numbers and some of their properties in Section
\ref{Sec:Pisot}.

We conclude with two notes.  First, the observations in the previous paragraph for affine IFSs on
$\mathbb{R}$ ($d=1$) apply \textit{mutatis mutandis} to affine IFSs (\ref{tauB}) in
$\mathbb{R}^d$, except for the results of Erd\H{o}s, which we extend in Sections \ref{Sec:2D} and \ref{Sec:higherdim} to $d\geq 2$. 
Second, a more general case of an affine IFS on $\mathbb{R}$ is prescribed by $\lambda\in(0,1)$, a subset $\{b_0,b_1\}\subset\mathbb{R}$, and the two affine maps
\begin{equation}
\tau_0(x) = \lambda(x + b_0)\quad \textrm{ and }\quad \tau_1(x) = \lambda(x+ b_1).
\end{equation}\label{generaltwotau}

\noindent We will return to IFSs of the form (\ref{generaltwotau}) in Sections \ref{Sec:Pisot} and \ref{Sec:InducedMeasures}.

\subsection{The general case}\label{Subsec:IFSGeneralCase}

An affine IFS in $\mathbb{R}^d$ is specified by a given $d \times d$
invertible and expansive matrix $A$ (the eigenvalues of $A$ are larger
than $1$ in absolute value) and a finite subset $B$ in
$\mathbb{R}^d$. The affine mappings are then indexed by the set $B$;
specifically
\begin{equation}
\tau_b(x) := A^{-1} (x + b),
\end{equation}\label{tauB}%
where $x \in \mathbb{R}^d$ and $b$ ranges over the set $B$. By
\cite{Hut81}, each such system (\ref{tauB}) has a natural equilibrium measure
$\mu_{(A, B)}$ with compact support $X_{(A,
B)}\subset\mathbb{R}^d$. 

               In what follows, the term \textit{rational} refers to the
following two restrictions (see \cite{DuJo07}):
\begin{enumerate}[(a)]
\item  The matrix $A$ has integral entries
\item  $B$ is a subset of $\mathbb{Z}^d$.    
\end{enumerate}
In the rational case, the fractal nature of the measure
$\mu_{(A, B)}$ is well understood, and there are recent papers, see
e.g., \cite{DuJo07} and  \cite{JoPe98}, which explore
further conditions on the pair $(A, B)$ which guarantee that
$L^2(\mu_{(A, B)})$ has an orthogonal basis of complex
exponentials. By this we mean that there is a subset $\Lambda$ of
$\mathbb{R}^d$ such that the set of complex exponential functions
\begin{equation*}
\{e_{\ell}(x) := \exp(2\pi i \ell\cdot x)\:|\: x \in X_{(A, B)} \text{ and } \ell \in \Lambda\}
\end{equation*} 
forms an orthogonal basis in $L^2(\mu_{(A, B)})$; in other words,
$L^2(\mu_{(A, B)})$ has a Fourier basis.  When such a Fourier basis
exists, we say that the index set $\Lambda$ is a \textit{spectrum} and $\mu_{(A,B)}$ is a \textit{spectral measure}.   

Even more generally, we know from \cite{Hut81} that for each
distribution $p=\{p_b\}_{b\in B}$ such that $p_b > 0$ and
$\sum_{b\in B} p_b = 1$,
there is a unique probability measure $\mu_{(A,B,p)}$ satisfying
\begin{equation}
\mu_{(A,B,p)} = \sum_{b\in B} p_b\: \mu_{(A,B,p)}\circ \tau_b^{-1},
\end{equation}\label{Eqn:MuABp}%
with
\begin{equation*}
\textrm{supp}(\mu_{(A,B,p)}) = X_{(A, B)}.
\end{equation*}
A direct calculation yields the following formula for the Fourier transform, which is valid for all $\xi\in\mathbb{R}^d$:
\begin{equation}
\hat{\mu}_{(A,B,p)}(\xi)
= 
\int e^{i2\pi \xi\cdot x}d\mu_{(A,B,p)}(x).
\end{equation}\label{Eqn:GeneralFTPart1}%
We use  $\int f(x)d(\mu\circ\tau^{-1})(x) = \int f(\tau(x))d\mu(x)$ to rewrite Equation (\ref{Eqn:GeneralFTPart1}):
\begin{equation}
\hat{\mu}_{(A,B,p)}(\xi)
=
\hat{\mu}_{(A,B,p)}((A^t)^{-1}\xi)\sum_{b\in B} p_b e^{i2\pi A^{-1}b\cdot \xi}.
\end{equation}\label{Eqn:GeneralFTPart2}%
If we define 
\begin{equation}
m_B(\xi) = \sum_{b\in B}p_b e^{i 2\pi  b\cdot \xi},
\end{equation}\label{Eqn:GeneralmB}%
and iterate the calculation in Equation (\ref{Eqn:GeneralFTPart2}), 
we obtain an infinite product formula for the Fourier transform, which converges because $m_B((A^t)^{-k}\xi)\rightarrow 1$ as $k\rightarrow\infty$:
\begin{equation}
\hat{\mu}_{(A,B,p)}(\xi) = \prod_{k=1}^{\infty} m_B((A^t)^{-k}\xi).
\end{equation}\label{Eqn:GeneralFTProduct}%
The goal of Section \ref{Sec:InducedMeasures} is to extend the infinite product formula (\ref{Eqn:GeneralFTProduct}) to a more general setting.

\begin{example}The simplest and best-known fractal, the middle-thirds Cantor set, is the atttractor for the IFS (\ref{generaltwotau}) with $\lambda = 1/3$, $b_0 = -1$, and $b_1 = 1$.  When $p_1 = p_2 = 1/2$,  Equation (\ref{Eqn:MuABp}) can be written
\begin{equation*}
\int \phi(x)\:d\mu(x) 
= \frac{1}{2}\Biggl(\int\phi\Bigl(\frac{x-1}{3}\Bigr)\:d\mu(x) 
+ \int\phi\Bigl(\frac{x+1}{3}\Bigr)\:d\mu(x) \Biggr)
\end{equation*}
for any $\phi$ which is integrable with respect to $\mu$.
In this case, Equation (\ref{Eqn:GeneralFTProduct}) becomes
\begin{equation*}
\hat{\mu}(\xi) = \prod_{n=1}^{\infty} \cos\Bigl(\frac{2\pi \xi}{3^n}\Bigr).
\end{equation*}
We have already mentioned that the result of Erd\H{o}s \cite{Erd39} implies that if $\lambda$ is the inverse of a Pisot number, then the Fourier transform of the associated Hutchinson measure does not tend to $0$ at infinity.  This result is a consequence of a property of Pisot numbers, which we will explain in Section \ref{Sec:Pisot}.  Here, our focus is on $\lambda = 1/3$;  we note that even though $1/3$ is not a Pisot number, the same type of result is even easier to obtain for $\hat{\mu}$: for any $m\in\mathbb{N}$,
\begin{equation*}
|\hat{\mu}(1)| = |\hat{\mu}(3^m)| > 0.  
\end{equation*}
As a result, $\mu$ is not absolutely continuous with respect to Lebesgue measure; in fact, $\mu$ is purely singular and continuous, as we will see in Section \ref{Subsec:Gradation}, Theorem \ref{Thm:PurelySingularORAC}.
\end{example}\label{Ex:Cantor1/3Product}%

The following proposition is \cite[Theorem 3.1 and Example 3.2]{DuJo07}:
\begin{proposition}
If $A$, $B$, and $p$ are as above, and if 
\begin{equation*}
\sum_{b\in B} \Bigl(p_b - \frac{1}{\# B}\Bigr)^2 > 0,
\end{equation*}
and $\sum_{b\in B}p_b = 1$, then the Hilbert space $L^2(\mu_{(A,B,p)})$ does \textbf{not} have an orthonormal basis (ONB) of complex exponentials  $e_{\ell}(x) = \exp(2\pi i \ell\cdot x)$ for any choice of points $\ell\in\mathbb{R}^d$.
\end{proposition}\label{Prop:DuJo07Theorem}%
\noindent As a result of Proposition \ref{Prop:DuJo07Theorem}, we will restrict to the equidistribution, i.e., $p_b = \frac{1}{\#B}$ for each $b\in B$, and we will use $\mu_{S}$ to denote the Hutchinson measure associated with the equidistribution $p$.

\subsection{Context of this paper}\label{Subsec:FocusContext}

       As we outlined in the Introduction, it is of great interest to
explore inverse spectral questions for a wider class of measures $\mu$
which arise in the theory of IFSs.  By ``inverse spectral questions'' we mean studying properties of an IFS determined by $S = \{\tau_i\}$ 
via the Fourier transform $\hat{\mu}_S$ or $\hat{\mu}_{(A,B)}$ in Equation  (\ref{Eqn:GeneralFTProduct}). 

The papers \cite{DuJo07} and \cite{JKS07a, JKS07b} study such questions for
various classes of affine IFSs.   In \cite{JKS07b}, we considered
orthogonal complex exponentials in $L^2(X, \mu_{(\lambda, \{\pm1\})})$ for the non-rational case in $\mathbb{R}$.  As we
noted there, relatively little is known about either the spectral or continuity properties of $\mu_{(A,B)}$.  Therefore, for the sake of simplicity in higher dimensions, we first consider systems (\ref{tauB}) where $A = \lambda^{-1}I$ and $B$ consists of standard basis vectors in $\mathbb{R}^d$ along with the zero vector.  

We introduce a 2D analogue of \cite{Erd39} in Section \ref{Sec:2D} and show that if $\lambda^{-1}$ is a Pisot number, then the 2D version of $\mu_{\lambda}$ is purely singular.  Our result in 2D is perhaps surprising when $\lambda > 2/3$
because, by \cite{JKS07a}, that is when our affine IFS does not have
gaps; i.e., the support of $\mu_{\lambda}$ is the planar closed set bounded by a
triangle with side lengths depending on $\lambda$.  In other words, the attractor of the IFS does not have a fractal structure.  In this case, we might expect $\mu_{\lambda}$ to be absolutely continuous with respect to Lebesgue measure, but for our special values of $\lambda$, $\mu_{\lambda}$ is not (Theorem \ref{Thm:Erdos2D}).  We can compare this result with Solomyak's result for 1D, which states that for almost every $\lambda\in (1/2, 1)$, $\mu_{\lambda}$ is absolutely continuous with respect to Lebesgue measure \cite{Sol95}.  In 1D, when $\lambda\in (1/2, 1)$, the fractal structure of the attractor also disappears---the attractor is a closed bounded interval. 

In 1D, the only known examples for which $\mu_{\lambda}$ is a spectral measure occur when $\lambda < 1/2$.  When $\lambda > 1/2$, the question of whether $\mu_{\lambda}$ is a spectral measure is much more difficult.  There are no known examples for which $\mu_{\lambda}$ is a spectral measure for $\lambda > 1/2$ in 1D, although there are many examples for which there are infinitely many orthonormal complex exponentials \cite[Theorem 2]{JKS07b}.   In higher dimensions, the spectral question is also difficult.  In \cite{DuJo07} it is conjectured that the only measures for which $\mu_{\lambda}$ is a spectral measure is the rational case (see Subsection \ref{Subsec:IFSGeneralCase}).  As a result, if the conjecture is true, the measures in Theorems \ref{Thm:Erdos2D} and \ref{Thm:ErdosdD} and the corresponding induced measures in Section \ref{Sec:InducedMeasures} would not be spectral measures. 

We note that the papers \cite{Erd39} and \cite{Sol95} are motivated by a study of the ``random'' geometric series $\sum\pm\lambda^n$ \cite{Kah85}.  The series is random because the plus/minus signs in front of $\lambda^n$ are the outcome of a sequence
of independent Bernoulli trials with two outcomes, such as flipping a coin.
That is, the assignment of the plus and minus signs is a random variable. In fact, for a fixed value of $\lambda$, the measure $\mu_{\lambda}$ is the density of the distribution of this random variable represented by the random power series.  We will return to this theme in Section \ref{Sec:InducedMeasures}.

\subsection{Helpful lemmas}\label{Subsec:HelpfulLemmas}

In the next three sections, we refer to the following lemmas about infinite products several times, so we state and prove them here for the sake of continuity and brevity in later proofs.

In $\mathbb{R}^d$, we work with the affine IFS associated with the set $B_d$, where
\begin{equation*}
B_d:=\{\textbf{0}, \textbf{e}_1, \ldots, \textbf{e}_d\},
\end{equation*}
and $\textbf{e}_i$ is the $i$th standard basis vector in $\mathbb{R}^d$.  
When $A = \lambda^{-1}I$, the Hutchinson measure $\mu_{\lambda}$ has Fourier transform
\begin{equation*}
\hat{\mu}_{\lambda}([\xi_1,\ldots, \xi_d]^t) 
= \prod_{n=1}^{\infty}\Biggl( \frac{1}{d+1}\Bigl(1 + \sum_{j = 1}^{d} e^{2\pi i \lambda^n  \xi_j}\Bigr)\Biggr),
\end{equation*}
a special case of Equation (\ref{Eqn:GeneralFTProduct}).

We restrict $[\xi_1, \ldots, \xi_d]^t$ to the line $L:=\{ \xi[1, \ldots, 1]^t:\xi\in\mathbb{R}\}$, and when we restrict $\hat{\mu}_{\lambda}$ to the line $L$, we write
\begin{equation*}
\begin{split}
& \hat{\mu}_{\lambda}([\xi,\ldots, \xi]^t)
 :=\hat{\mu}_{\lambda, B, W}(\xi) \\
& = \prod_{n=1}^{\infty} \Biggl(\frac{1}{d+1}\Bigl(1 +d e^{2\pi i \lambda^n \xi}\Bigr)\Biggr)
= \prod_{n=1}^{\infty}  m_{B,W}(\lambda^n \xi)
\end{split}
\end{equation*}
We will often work with estimates involving only the real part of $m_{B,W}$, which accounts for the following lemma about the cosine function.
\begin{lemma} Let $d\in\mathbb{N}$ and $\lambda \in (0,1)$.  There exists $N\in\mathbb{N}$ and $C >0$ such that 
\begin{equation*}
\prod_{n=N}^{\infty}  \Biggl(\frac{1}{d+1}\Bigl(1 +d \cos(2\pi  \lambda^n)\Bigr)\Biggr) > C.
\end{equation*}
\end{lemma}\label{Lemma:Taylord}%
\textit{Proof:  }We start by taking logarithms to convert the product to a sum: 
\begin{equation*}
\prod_{n = N}^{\infty} \Biggl( \frac{1}{d+1}\Bigl( 1 + d\cos(\lambda^n)\Bigr) \Biggr) > C 
\end{equation*}
if and only if
\begin{equation*}
\sum_{n=N}^{\infty} \ln\Biggl( \frac{1}{d+1}\Bigl( 1 + d\cos(\lambda^n)\Bigr) \Biggr) > \ln(C). \end{equation*}
We will use the Taylor series expansion of $1 + d\cos(x)$ around $x=0$
and the Taylor series of $\ln(x)$ around $x=1$ to find $C$.  First,
\begin{equation*}
\frac{1}{d+1} + \frac{d}{d+1}\cos(\lambda^n) = 1 - \sum_{k = 1}^{\infty} \frac{(-1)^{k+1}d\lambda^{2kn}}{(d+1)(2k)!}.
\end{equation*}
Let 
\begin{equation*}\varepsilon_n = \sum_{k = 1}^{\infty} \frac{(-1)^{k+1}d\lambda^{2kn}}{(d+1)(2k)!};\end{equation*}
we note that $\varepsilon_n \geq 0$ for all $n$.  Choose $N_1\in\mathbb{N}$ such that for all $n > N_1$, 
$\varepsilon_n < \frac{\lambda^2}{2}$.
Next, choose $N_2 > N_1$ such that for all $n > N_2$, 
$\frac{\lambda^2}{2} < 1$. 
Now examine the Taylor series expansion for $\ln(\frac{1}{d+1} + \frac{d}{d+1}\cos(\lambda^n))$, which is valid for all $n > N_2$:
\begin{equation*} \ln\Bigl(\frac{1}{d+1} + \frac{d}{d+1}\cos(\lambda^n)\Bigr) = \ln(1 - \varepsilon_n) = - \sum_{k = 1}^{\infty}\frac{(\varepsilon_n)^k}{k}.
\end{equation*}
Finally, choose $N > N_2$ such that for all $n > N$, 
\begin{equation*} \sum_{k = 1}^{\infty}\frac{(\varepsilon_n)^k}{k}   \leq 2\varepsilon_n.\end{equation*}
Now,
\begin{equation*}
\sum_{n = N}^{\infty} \ln(1 - \varepsilon_n)
= \sum_{n = N}^{\infty} \Biggl(-\sum_{k = 1}^{\infty}\frac{(\varepsilon_n)^k}{k}\Biggr) 
\geq \sum_{n = N}^{\infty} -2\varepsilon_n \geq \sum_{n = N}^{\infty} -\lambda^{2n} 
= \ln(C).
\end{equation*}
$ $\hfill$\Box$

We note that Lemma \ref{Lemma:Taylord} follows the same lines as \cite[Lemma 2]{JKS07b}.  Lemma \ref{Lemma:Taylord} tells us that the only way that the infinite product 
\[ \prod_{n=1}^{\infty}  \Biggl(\frac{1}{d+1}\Bigl(1 +d \cos(2\pi  \lambda^n)\Bigr)\Biggr)\]
can be $0$ is if one of its
factors is $0$.

\begin{lemma}
Let $\lambda\in (0,1)$.  In $\mathbb{R}^d$, $d > 1$, the product
\begin{equation}
\prod_{n=1}^{\infty}|m_{B,W}(\lambda^{n})|^2
\end{equation}\label{Exp:ProductForLemmamBLlambda}%
is nonzero.
\end{lemma}\label{Lemma:mBLlambda}%
\textit{Proof:  }
Choose $N$ as in Lemma \ref{Lemma:Taylord} so that
\begin{equation*}
\prod_{n=N}^{\infty}\Bigl(\frac{1}{d+1}\Bigr)^2\Bigl(1 + d\cos(2\pi \lambda^n)\Bigr)^2 > 0,
\end{equation*}
and split the infinite product (\ref{Exp:ProductForLemmamBLlambda}) into two parts:
\begin{equation*}
\prod_{n=1}^{\infty}|m_{B,W}(\lambda^{n})|^2
=
\prod_{n=1}^{N-1}|m_{B,W}(\lambda^{n})|^2
\prod_{n=N}^{\infty}|m_{B,W}(\lambda^{n})|^2.
\end{equation*}
The finite product
\begin{equation*}
\prod_{n=1}^{N-1}|m_{B,W}(\lambda^n)|^2
\end{equation*}
is a positive constant because $m_{B,W}(x) = \frac{1}{d+1}\Bigl( 1 + d e^{2\pi i x}\Bigr)$, and
\begin{equation*}
1 + de^{2\pi i  x} = 0\textrm{ if and only if } e^{2\pi i x} = -\frac{1}{d},
\end{equation*}
which is impossible when $d > 1$.
Now
\begin{equation}
\begin{split}
&\prod_{n=N}^{\infty}|m_{B,W}(\lambda^{n})|^2
\geq
\prod_{n=N}^{\infty}[\textrm{Re}(m_{B,W}(\lambda^n))]^2\\ 
& = \prod_{n=N}^{\infty}\Bigl(\frac{1}{d+1}\Bigr)^2\Bigl(1 + d\cos(2\pi \lambda^n)\Bigr)^2 
> 0.
\end{split}
\end{equation}\label{replacewithcosine}%

$ $\hfill$\Box$

\begin{lemma}\label{Lemma:CountableTheta}
There are countably many values of $\theta$ such that
\[\prod_{n=0}^{\infty} \Bigl(\frac{1}{d+1}\Bigr)^2\Bigl(1 + d \cos(2\pi \theta^n) \Bigr)^2=0.\]
\end{lemma}

\section{Pisot numbers and Hutchinson measures}\label{Sec:Pisot}

With view to later use, we will explain the techniques and results of  Erd\H{o}s's 1939 paper \cite{Erd39} in this section.  We change his notation to fit our own in subsequent sections.  Erd\H{o}s's  techniques will be used in our later results.  
\subsection{Elementary results about Pisot numbers}

Let $\alpha$ be an algebraic integer---that is, $\alpha$ is the root
of a polynomial whose leading coefficient is $1$ and the rest of the
coefficients are all integers.  The algebraic integer $\alpha$ is a
\textit{Pisot number} (also known as a Pisot-Vijayaraghavan number or
PV-number) if $\alpha > 1$ and all the Galois conjugates of $\alpha$
have modulus less than $1$.  We will denote the Galois conjugates of
$\alpha=\alpha_1$ as follows: $\{\alpha_2, \ldots, \alpha_n\}$.

Large powers of Pisot numbers are ``\textit{almost}'' integers: consider the
expression $\sum_{i = 1}^{n} \alpha_i^k$, which is an integer by Lemma
\ref{Lemma:Pisot}.  As $k\rightarrow\infty$, all the terms in the sum
except the first tend to $0$.  Therefore, for large $k$,
\begin{equation*} 
\sum_{i = 1}^{n} \alpha_{i}^{k} \approx \alpha_1^k.
\end{equation*}
The expression $\sum_{i = 1}^{n} \alpha_i^k$ is sometimes called the
\textit{trace} or the \textit{spur} of $\alpha^k$
\cite{Ca57}. 

\begin{lemma}
Suppose $\alpha$ is an algebraic integer.  
For every $k\in\mathbb{Z}^{+}$, $\sum_{i = 1}^{n} \alpha_i^k$ is an
integer.\end{lemma}\label{Lemma:Pisot}%

\noindent For the convenience of the reader, we present two proofs of
Lemma \ref{Lemma:Pisot}.  While these proofs are known, the ideas from the discussion below are relevant for our later considerations.

\noindent\textit{Proof \#1 of Lemma \ref{Lemma:Pisot}.  }   This proof relies on the Fundamental Theorem of Symmetric Polynomials.  We will use the statement of this theorem from \cite{PoDi}:
\begin{theorem}\rm \cite[p. 36]{PoDi}
\it Every symmetric polynomial in $x_1, \ldots, x_n$ over a field $F$ can be written as a polynomial over $F$ in the elementary symmetric functions $\sigma_1, \ldots, \sigma_n$.  If the coefficients of the first are rational integers, the same is true of the second.
\end{theorem}\label{Thm:SymPol}%
The $i$th elementary symmetric function $\sigma_i$ is the ``sum of all products of $i$ different $x_j$'' \cite{PoDi}.  For example, 
\begin{equation*}
\sigma_1 = x_1 + \cdots + x_n
\end{equation*} and 
\begin{equation*}
\sigma_2 = x_1x_2 + x_1x_3 + \cdots + x_2x_3 + \cdots + x_{n-1}x_n.
\end{equation*}
Now let $p(x)\in\mathbb{Z}[x]$ be the minimal polynomial of $\alpha$:
\begin{equation*}
p(x) = x^n + a_1 x^{n-1} + \cdots + a_{n-1}x + a_n = \prod_{i = 1}^n(x-\alpha_i).
\end{equation*}
The coefficient $a_i$ is precisely $(-1)^i$ multiplied by the $i$th
elementary symmetric function in $\alpha_1, \ldots, \alpha_n$
\cite[p. 39]{PoDi}.  Now note that $\sum_{i = 1}^{n} \alpha_i^k$ is a
symmetric polynomial in $\alpha_1, \ldots, \alpha_n$ with coefficients
in $\mathbb{Z}$.  But by Theorem \ref{Thm:SymPol}, we can write
$\sum_{i = 1}^{n} \alpha_i^k$ as a linear combination of the
elementary symmetric functions in $\alpha_1, \ldots, \alpha_n$ with
coefficients in $\mathbb{Z}$.  Since the elementary symmetric
functions in $\alpha_1, \ldots, \alpha_n$ are the integers $a_1,
\ldots, a_n$, the sum $\sum_{i=1}^{n}\alpha_i^k$ must also be an
integer.  \hfill$\Box$

\noindent\textit{Proof \#2 of Lemma \ref{Lemma:Pisot}.  } Fix
$k\in\mathbb{N}$, and let $L = \mathbb{Q}(\alpha_1, \ldots, \alpha_n)$
be the splitting field of the minimal polynomial of $\alpha$.  Thus,
$L$ is a Galois extension since $\mathbb{Q}$ has characteristic $0$.
If $\sigma$ is \textit{any} element of $\textrm{Gal}(L, \mathbb{Q})$,
then $\sigma$ is a field automorphism, and $\sigma$ permutes the roots
of the minimal polynomial of $\alpha$.  Therefore
\begin{equation*} 
\sigma(\sum_{i=1}^n \alpha_i^k) 
= \sum_{i = 1}^n (\sigma(\alpha_i))^k 
= \sum_{i=1}^n \alpha_i^k,
\end{equation*} 
and $\sum_{i=1}^n \alpha_i^k$ is invariant under the action of the
entire Galois group $\textrm{Gal}(L, \mathbb{Q})$.  By the Galois
correspondence, the entire Galois group fixes \textit{only} the base
field $\mathbb{Q}$.  Therefore $\sum_{i=1}^n \alpha_i^k$ is an
algebraic integer which is also a rational integer---that is, the sum
is in $\mathbb{Z}$. \hfill$\Box$

$ $\hfill$\Box\Box$

With Lemma \ref{Lemma:Pisot} in hand, we can now see that successively
higher powers of $\alpha$ get closer and closer to integers.  In fact,
we can use the result in Lemma \ref{Lemma:Pisot} to control how close
$\alpha^k$ is to an integer.  In what follows, we will consider
$\alpha^k$ to belong to $\mathbb{R}/\mathbb{Z}\cong[0, 1)$, so
that integers are represented by $0$.

\begin{lemma}
There exists $\theta\in(0,1)$ such that for all $k\in\mathbb{N}$
the distance mod $1$ between $\alpha^k$ and $0$ is less than
$\theta^k$.
\end{lemma}\label{Lemma:GeometricBound}%
\textit{Proof:  }
Since $\sum_{i=1}^n \alpha_i^k\in\mathbb{Z}$ for any $k\in\mathbb{N}$, we have
\begin{equation*}
\begin{split}
\textrm{dist}(\alpha^k, \mathbb{Z})
& = \textrm{dist}(\alpha^k, 0)  = \Big|\sum_{i = 1}^{n} \alpha_i^k -
\alpha^k\Big| 
= \Big|\sum_{i = 2}^{n}\alpha_i^k\Big|\\ 
& \leq \sum_{i=2}^n
|\alpha_i|^k \leq (n-1)\max_{2\leq i \leq n} |\alpha_i|^k.
\end{split}
\end{equation*}
Suppose $\max_{2\leq i \leq n} |\alpha_i| = |\alpha_2|$.  We work with numbers and not their equivalence classes mod $1$ first.  Since $|\alpha_2|< 1$, we can choose $N\in\mathbb{N}$ and $\theta_1\in (0,1)$ such that
\begin{equation*}
(n-1)|\alpha_2 |^N < [(\theta_1)^{1/N}]^N.
\end{equation*}
In particular, we know that $|\alpha_2| < (\theta_1)^{1/N}$, so for all $j \geq 0$,
\begin{equation*}
(n-1)|\alpha_2 |^{N+j}  < [(\theta_1)^{1/N}]^{N+j}.
\end{equation*}
Now consider the equivalence classes (which we denote by $\overline{*}$ since all numbers are real)
\begin{equation*}
\{ 
\overline{(n-1)|\alpha_2|}, 
\overline{(n-1)|\alpha_2|^2},
\ldots,
\overline{(n-1)|\alpha_2|^{N-1}}
\};
\end{equation*}
choose $\theta_2 \in (0,1)$ so that
\begin{equation*}
\max_{1\leq i \leq N-1} \overline{(n-1) |\alpha_2|^j} < \theta_2^{N-1}.
\end{equation*}
Then for each $j\in \{1, 2, \ldots, N-1\}$,
\begin{equation*}
\overline{(n-1) |\alpha_2|^j} < \theta_2^{N-1} \leq \theta_2^j.
\end{equation*}
Set $\theta = \max\{(\theta_1)^{1/N}, \theta_2\}$; this value of $\theta$ satisfies the conclusion of the lemma.
\hfill$\Box$

\subsection{Erd\H{o}s's proof on Bernoulli convolution measures}\label{Subsec:ErdosProof}

Erd\H{o}s's proof of his 1939 theorem on Bernoulli convolution
measures is short and elegant.  In the paper, he assumes Lemma
\ref{Lemma:Pisot}, Theorem \ref{Thm:SymPol}, and
Lemma \ref{Lemma:GeometricBound} are familiar to the reader; the heart of
his proof consists of three lines of inequalities.  We state the theorem in terms of affine IFSs; the Hutchinson measure in this case is often called the infinite Bernoulli convolution measure.

\begin{theorem}\label{Thm:Erdos}\rm (Erd\H{o}s \cite{Erd39})\it\:\:
Suppose $\alpha$ is a Pisot number and $\lambda = \alpha^{-1}$.  Let $\mu_{\lambda}$ denote the Hutchinson measure for the IFS (\ref{generaltwotau}).  The Fourier transform $\hat{\mu}_{\lambda}$ is given by (\ref{Eqn:GeneralFTProduct}) with $A = \lambda^{-1}$, $B = \{-1, 1\}$, and $p = \{1/2, 1/2\}$.  

Then there exists $C > 0$ such that $|\hat{\mu}_{\lambda}(\alpha^k)| > C$ for all $k\in\mathbb{N}_0$.
\end{theorem}
\begin{remark}In Example \ref{Ex:Cantor1/3Product}, we saw that when $\alpha = 3$ and $\lambda = 1/3$, $|\hat{\mu}_{1/3}(3^n)| = |\hat{\mu}_{1/3}(1)| > 0$ for all $n$.
\end{remark}

\textit{Proof:  }
We use our own notation which is similar to but not identical to Erd\H{o}s's.  The Fourier transform of the
measure $\mu_{\lambda}$ arising from the IFS in (\ref{generaltwotau}) with $b_0 = -1, b_1 = 1$ is
\begin{equation}
\hat{\mu}_{\lambda}(\xi) = \prod_{n=1}^{\infty} \cos(2\pi \lambda^n \xi).
\end{equation}\label{Eqn:BernoulliConvolution}%
Fix $k\in\mathbb{N}$.  In what follows, we will exploit the fact that $\alpha^k$ is ``almost'' an integer:
\begin{equation}
\begin{split}
\hat{\mu}_{\lambda}(\alpha^k) 
& = \prod_{n=1}^{\infty} \cos(2\pi \lambda^n \alpha^k)\\
& = \prod_{n=1}^{\infty} \cos(2\pi \lambda^{n-k}).
\end{split}
\end{equation}\label{Eqn:OneDProductPart1}%
We now split the single product into two products, one with positive powers of $\lambda$ and one with positive powers of $\alpha$:
\begin{equation}
\begin{split}
\hat{\mu}_{\lambda}(\alpha^k) 
& = \prod_{n=1}^{\infty} \cos(2\pi \lambda^n)\prod_{n = -k+1}^0\cos(2\pi \lambda^{n})\\
& = \prod_{n=1}^{\infty} \cos(2\pi \lambda^n)\prod_{n = 0}^{k-1}\cos(2\pi \alpha^{n}).
\end{split}
\end{equation}\label{Eqn:OneDProductPart2}%
By Lemma \ref{Lemma:Taylord}, the infinite product in Equation (\ref{Eqn:OneDProductPart2}) is a positive constant.  

By Lemma \ref{Lemma:GeometricBound} we can choose $\theta\in(0,1)\backslash\{1/4\}$ and $N\in\mathbb{N}$ such that
\begin{equation*}
\alpha^n (\text{ mod } 1) < \theta^n < 1/4 \text{ for all }n \geq N.
\end{equation*}
We make sure $\theta\neq 1/4$ so that the product $\prod_{n=0}^{\infty} \cos(2\pi \theta^n)$ is nonzero, and we require that $\theta^n < 1/4$ because $\cos(2\pi x)$ is a decreasing function of $x$ in $(0, 1/4)$.  

\bigskip

\noindent\textbf{Case 1:  }$k -1 \geq N$

\bigskip

Set
\begin{equation}
C = \prod_{n=0}^{N-1} |\cos(2\pi \alpha^{n})|;
\end{equation}\label{Eqn:ErdosC}%
$C$ is nonzero by Lemma \ref{Lemma:FiniteProductNonzeroErdos} below.

We now split the finite product in Equation (\ref{Eqn:OneDProductPart2}) into two products and use $\theta$ to create a lower bound:
\begin{equation}
\begin{split}
\prod_{n = 0}^{k-1}|\cos(2\pi \alpha^{n})| 
= C \prod_{n = N}^{k-1}|\cos(2\pi \alpha^{n})|     
> C \prod_{n = N}^{k-1}|\cos(2\pi \theta^{n})|.
\end{split}
\end{equation}\label{Eqn:OneDBounds}%
Since every term of the product $\prod_{n = 0}^{\infty}|\cos(2\pi \theta^{n})|$ lies between $0$ and $1$, we can bound the last finite product (which depends on $k$) below with an infinite product (which does not depend on $k$):
\begin{equation*}
 \prod_{n = N}^{k-1}|\cos(2\pi \theta^{n})| \geq \prod_{n = 0}^{\infty}|\cos(2\pi \theta^{n})|, 
\end{equation*}
and now we have a positive lower bound for the sequence 
$\{|
\hat{\mu}_{\lambda}(\alpha^k)|
\}_{k = N}^{\infty}$.  

\bigskip

\noindent\textbf{Case 2:  }$0 \leq  k \leq N$

\bigskip

Referring back to Equation (\ref{Eqn:OneDProductPart2}), we know that the infinite product is positive by Lemma \ref{Lemma:Taylord}, and we know that the finite product is positive by Lemma \ref{Lemma:FiniteProductNonzeroErdos}.

\bigskip

Combining Cases 1 and 2, we obtain a postive lower bound for the sequence $\{|\hat{\mu}_{\lambda}(\alpha^k)|\}_{k = 0}^{\infty}$.
\hfill$\Box$

\begin{lemma}
Let $\alpha$ be a Pisot number.  Then for any $j\in\mathbb{N}$,
\begin{equation}
\prod_{n=0}^j |\cos(2\pi \alpha^n)| > 0.
\end{equation}\label{Ineq:FiniteProductPositiveErdos}%

\end{lemma}\label{Lemma:FiniteProductNonzeroErdos}%

\textit{Proof:  }
Since $\alpha$ is a Pisot number, $\alpha^n$ is not rational for any $n\in\mathbb{N}$.  The only zeros of $\cos(2\pi x)$ are rational.
\hfill$\Box$

\begin{remark}The techniques of Erd\H{o}s's theorem will carry over to our theorems for affine IFSs in $\mathbb{R}^d$, $d > 1$.  In particular, inequalities very much like (\ref{Ineq:FiniteProductPositiveErdos}) in Lemma \ref{Lemma:FiniteProductNonzeroErdos} will return again in modified forms in 
\begin{itemize}
\item Inequality (\ref{Ineq:FiniteProductNonzero2D}) in Lemma \ref{Lemma:FiniteProductNonzero2D} for Theorem \ref{Thm:Erdos2D} in Section \ref{Subsec:FixedDirection2D}
\item Inequality (4.21) in Theorem \ref{Thm:GenDir2D}, Section \ref{Subsec:GeneralDirection2D}
\item Inequality (2.26) in Lemma \ref{Lemma:mBLlambda} for Theorem \ref{Thm:ErdosdD} in Section \ref{Sec:higherdim}.
\end{itemize}  
The common theme here is that we work with finite products with
factors indexed from $0$ to $j$, where $j\in\mathbb{N}$.  In all cases, the factors in the products are less than $1$, and in all cases, the conclusion is that for
each $j$, there is a positive lower bound for the product which may be chosen independently of $j$. In fact, the products
converge to a fixed and positive lower bound as $j\rightarrow\infty$.
\end{remark}

\begin{remark}\textbf{Extensions of Theorem \ref{Thm:Erdos} in 
$\mathbb{R}$.  }

\noindent While part of our focus is extending Erd\H{o}s's theorem to higher dimensions, we note that the theorem has been extended to other affine IFSs on $\mathbb{R}$ in \cite{LNR01}.  In \cite[Theorem 5.1]{LNR01}, the authors extend Erd\H{o}s's theorem to affine IFSs on $\mathbb{R}$ in the following way.  If $\{b_i\}\subset\mathbb{Q}$, $\{p_i\}$ is any set of positive probability weights, and $\lambda$ is the reciprocal of a Pisot number, then the invariant measure of the affine IFS defined by $\{p_i\}$ and $\{ \lambda x + b_i\}_{i = 1}^N$ is singular. 
\end{remark}\label{Rem:RealExtensions}%

\subsection{The dichotomy in the class of Hutchinson measures}\label{Subsec:Dichotomy}

Recall that a measure $\mu$ on $\mathbb{R}^d$ has a canonical decomposition into the sum of three parts: 
\begin{itemize}
\item $\mu_{ac}$, the absolutely continuous part (with respect to Lebesgue measure on $\mathbb{R}^d$),
\item $\mu_{pp}$, the atomic part, and
\item $\mu_{cs}$, the purely continuous and singular part.
\end{itemize}

\begin{definition}\label{Defn:PurelySingular}
A measure $\mu$ on $\mathbb{R}^d$ is said to be purely continuous
and singular if in its decomposition the two components $\mu_{ac}$ and
$\mu_{pp}$ are zero.
\end{definition}

We note that by Wiener's theorem, Hutchinson measures arising from affine IFSs do not have any point masses \cite[Corollary 6.6]{DuJo06c}.  In fact, more is known:
\begin{theorem}If $\mu$ is a Hutchinson measure defined by (\ref{Eqn:MuABp}), and if $\mu = \mu_{ac} + \mu_{cs}$, then either $\mu_{ac}=0$ or $\mu_{cs} = 0$.
\end{theorem}\label{Thm:PurelySingularORAC}%

\begin{remark}In particular, if $\lambda$ is the reciprocal of a Pisot number and if $\mu = \mu_{\lambda}$ is the infinite Bernoulli convolution measure in Theorem \ref{Thm:Erdos}, $\mu_{\lambda}$ cannot have an absolutely continuous part by Theorem \ref{Thm:PurelySingularORAC}, so $\mu_{\lambda}$ is purely continuous and singular.  Erd\H{o}s says that $\mu_{\lambda}$ is ``purely singular'' in \cite{Erd39}.
\end{remark}

\begin{remark}
Theorem \ref{Thm:PurelySingularORAC} was known to Jessen and Wintner \cite[Theorem 11]{JeWi35} for Bernoulli convolution measures on $\mathbb{R}$; a proof of this theorem is outlined in \cite[Proposition 3.1, p. 42-43]{PSS98}.  We will prove the theorem as stated here in full for the benefit of the reader.
\end{remark}

Throughout, we consider only positive Borel measures in $\mathbb{R}^d$ which are all compactly supported and therefore $\sigma$-finite.  We begin by recalling the following definitions and terminology.
\begin{definition}\cite[p. 120]{Rud87}
\begin{enumerate}[(a)]
\item The measure $\mu$ is \textit{concentrated} on a set $A$ if for every Borel set $E\subset\mathbb{R}^d$, \[\mu(A\cap E) = \mu(E).\]
\item The notation $\mu_s\perp m$ means that there exists a pair of disjoint sets $X_{s}$ and $L$ such that $\mu_s$ is concentrated on $X_{s}$ and $m$ is concentrated on $L$.  In this case, we say that $\mu_{s}$ and $m$ are \textit{mutually singular}.  In the present context, this means that $\mu_{s}$ has no component which is absolutely continuous with respect to $m$.  We note that a measure can be concentrated on many different sets.
\item The \textit{support} of the measure $\mu$, denoted $\textrm{supp}(\mu)$, is the smallest closed set on which $\mu$ is concentrated.  An equivalent definition is that $\textrm{supp}(\mu) = A$ if and only if for all $\phi\in C_{c}(\mathbb{R}^d)$,
\begin{equation*}
\int_{\mathbb{R}^d} \phi d\mu = \int_{A} \phi d\mu.
\end{equation*}
When there exists $f\in L^1(m)$ such that $\mu = f \,dm$, then 
$\textrm{supp}(\mu) = \textrm{supp}(f)$, where the last support refers to that of a function.
\end{enumerate}
\end{definition}
A special case of the Lebesgue-Radon-Nikodym Theorem \cite[Theorem 6.10, p. 121]{Rud87} tells us that there exist \textit{unique} measures $\mu_{ac}$ and $\mu_{s}$ such that
\begin{equation}
\mu = \mu_{ac} + \mu_s,
\end{equation}\label{Eqn:Rudin6.10a}%
and if $m$ denotes Lebesgue measure on $\mathbb{R}^d$, then
\begin{equation}
\mu_{ac}    \ll  m \qquad\textrm{ and } \qquad \mu_{s} \perp \mu_{ac}.
\end{equation}\label{Eqn:Rudin6.10b}%

Let $(A, B)$ define our IFS as in (\ref{tauB}), where $\mathbf{0}\in B$.  Let $H$ denote the transformation on measures
\begin{equation}
H\mu := \sum_{b\in B} p_b \mu \circ\tau_b^{-1};
\end{equation}\label{Eqn:Hmu}%
our goal is to show that 
\[H\mu_{ac} = \mu_{ac} \textrm{ and }H\mu_s = \mu_s.\]  To show these two equations are true, 
we compute Radon-Nikodym derivatives.  We start with the more straightforward case---the Radon-Nikodym derivative of $H\mu_{ac}$.

\begin{lemma}
If $d\mu_{ac} = f(x) dm$ where $f\in L^{1}(m)$, $f\geq 0$, then the Radon-Nikodym derivative of $H\mu_{ac}$ is given by $g\in L^1(m)$, where
\begin{equation*}
g(x) =  |\det(A)|\sum_{b\in B}  p_b f\circ\tau_b^{-1}(x).
\end{equation*}
\end{lemma}\label{Lemma:RNDerHmuac}%
\textit{Proof:  }
Let $\phi$ be a test function:
\begin{equation*}
\begin{split}
\int \phi(x) dH\mu_{ac}(x) 
& = \int \phi(x) d\Bigl(\sum_{b\in B} p_b \mu_{ac}\circ\tau_{b}^{-1}(x)\Bigr)\\
& = \sum_{b\in B} p_b \int \phi(\tau_b(x))\:d\mu_{ac}(x)\\
& = \sum_{b\in B} p_b \int \phi(\tau_b(x))f(x)\:dx.
\end{split}
\end{equation*}
Now, change variables and let $\tau_b(x) = A^{-1}(x + b) = y$.  Then $dx = |\det(A)|dy$ and $x = \tau_b^{-1}(y) = A y - b$.  We rewrite the integral
\begin{equation*}
\int \phi(x) dH\mu_{ac}(x)  = \sum_{b\in B} p_b\int\phi(x) f(\tau_b^{-1}(x)) |\det(A)| dx,
\end{equation*}
and the result follows.
\hfill$\Box$

In order to compute $DH\mu_{s}$, the Radon-Nikodym derivative of $H\mu_s$, we use \cite[Theorem 7.14, p. 143]{Rud87}.  First,  we give a definition and a lemma about (\ref{tauB}).

\begin{definition}\cite[Definition 7.9, p. 140]{Rud87}
Let $\{E_i\}_{i=1}^{\infty}$ be a collection of subsets of $\mathbb{R}^d$.  We say that $\{E_i\}$ \textit{``shrinks nicely''} to $x$ if there exists $\alpha > 0$ and there exists a sequence $r_i \rightarrow 0$ such that
\begin{enumerate}
\item $E_i \subset B(x,r_i)$
\item$m(E_i) \geq \alpha \,m(B(x,r_i))$.
\end{enumerate}
The set $E_i$ does not have to contain the point $x$.
\end{definition}\label{Defn:SN}%

Set 
\begin{equation}
D\mu(x) = \lim_{i\rightarrow\infty}\frac{\mu(E_i)}{m(E_i)}
\end{equation}\label{Eqn:RNDer}%
when this limit exists.  Rudin shows that $D\mu$ agrees with the Radon-Nikodym derivative $\frac{d\mu}{dm}$ when the latter exists.  In particular, $D\mu(x)$ is independent of the choice of sets $\{E_i\}$ which shrink nicely to $x$; we use this fact in Lemma \ref{Lemma:HmuSSingular} below.

\begin{lemma}
Let $x\in\mathbb{R}^d$ and let $\{E_i\}$ shrink nicely to $x$.  Then $\{\tau_b^{-1}(E_i)\}$ shrinks nicely to $\tau_{b}^{-1}(x)$.
\end{lemma}\label{Lemma:TausShrinkNicely}%
\textit{Proof:  }We are given  (1) and (2) in Definition \ref{Defn:SN} above, where $r_i\rightarrow 0$.  If $x\in E_i \subset B(x,r_i)$, then $\tau_b^{-1}(x) \in \tau_b^{-1}(E_i)\subset \tau_b^{-1}(B(x,r_i))$.  We know that $A^{-1}$ is contractive and that all the eigenvalues of $A$ are greater than $1$ in absolute value.  Therefore for each $i$ there exist positive numbers $\rho_i, R_i, a, b$ such that
\begin{itemize}
\item $0 < \rho_i < R_i$ 
\item $R_i\rightarrow 0$ 
\item $B(\tau_b^{-1}(x), \rho_i) \subset \tau_b^{-1}(B(x,r_i)) \subset B(\tau_{b}^{-1}(x), R_i)$
\item $a \leq \frac{R_i}{\rho_i} \leq b$ for all $i$
\end{itemize}
Also, there exists $\delta > 0$ depending only on the spectrum of $A$ such that
\begin{equation*}
m(\tau_b^{-1}(E_i)) = |\det(A)| m(E_i) \geq \alpha |\det(A)| m(B(x,r_i)) \geq \alpha \delta m(B(\tau_b^{-1}x, R_i)).
\end{equation*}
Since $R_i\rightarrow 0$, the result follows.
\hfill$\Box$

\begin{theorem}\rm \cite[Theorem 7.14, p. 143]{Rud87}
\it
Suppose that to each $x\in \mathbb{R}^d$ is associated some sequence $\{E_i(x)\}$ which shrinks to $x$ nicely, and that $\mu$ is a complex Borel measure on $\mathbb{R}^d$.  Let $\mu = f\:dm + d\mu_s$ be the Lebesgue decomposition of $\mu$ with respect to $m$.  Then 
\begin{equation*}
\lim_{i\rightarrow\infty} 
\frac{\mu(E_i(x))}{m(E_i(x))} = f(x)\:\:a.e. [m].
\end{equation*}
In particular, $\mu \perp m$ if and only if 
\begin{equation*}D\mu(x) = 0\:\: a.e. [m].
\end{equation*}
\end{theorem}

\begin{lemma}
Let $H\mu_s$ be as in (\ref{Eqn:Rudin6.10a}) and (\ref{Eqn:Hmu}) above. Then
\begin{equation*}
DH\mu_s(x)= 0 \:\:a.e. [m].
\end{equation*}
\end{lemma}\label{Lemma:HmuSSingular}
\textit{Proof:  }
It is enough to show that $D(\mu_s\circ\tau_b^{-1}) = 0$ a.e. [$m$].  Suppose $\{E_i(x)\}$ shrinks nicely to $x$.  To show that the Radon-Nikodym derivative of $\mu_s\circ\tau_{b}^{-1}$ is $0$, we show that 
\begin{equation}
\lim_{i\rightarrow\infty} \frac{(\mu_s\circ\tau_{b}^{-1})(E_i)}{m(E_i)}=0.
\end{equation}\label{Eqn:OriginalLimitDHmus}%
However, we can choose $c_1, c_2\in \mathbb{R}^{+}$ such that 
\begin{equation*}
c_1 m(\tau_b^{-1}(E_i))
 \leq 
m(E_i) 
\leq 
c_2 m(\tau_b^{-1}(E_i)),
\end{equation*}
so Equation (\ref{Eqn:OriginalLimitDHmus}) is true if and only if
\begin{equation}
\lim_{i\rightarrow\infty} \frac{\mu_s(\tau_b^{-1}(E_i))}{m(\tau_{b}^{-1}(E_i))}=0.
\end{equation}\label{Eqn:LastLimitDHmus}%
But $\{\tau_b^{-1}(E_i)\}$ shrinks nicely to $\tau_b^{-1}(x)$ by Lemma \ref{Lemma:TausShrinkNicely}, and $\mu_s$ is singular with respect to Lebesgue measure $m$---that is, $D\mu_s = 0$ a.e. [$m$].  Therefore the limit in (\ref{Eqn:LastLimitDHmus}) is $0$.\hfill$\Box$

We know from Lemma \ref{Lemma:HmuSSingular} that $H\mu_s$ is singular with respect to Lebesgue measure, and we know from Lemma \ref{Lemma:RNDerHmuac} that $H\mu_{ac}$ is absolutely continuous with respect to Lebesgue measure.  One consequence is that $H\mu_{s} \perp H\mu_{ac}$. The theorem of Hutchinson \cite{Hut81} tells us that $\mu$ is the unique solution to $\mu = H\mu$, so
\begin{equation*}
\mu = \mu_{ac} + \mu_{s} = H\mu_{ac} + H\mu_s = H\mu.
\end{equation*}
 By the Lebesgue-Radon-Nikodym theorem, the decomposition of $\mu$ into singular and absolutely continuous parts is unique, so 
\begin{equation}
\mu_{ac} = H\mu_{ac}
\end{equation}\label{Eqn:MuACisHMuAC}%
and 
\begin{equation*}
\mu_{s} = H\mu_{s}.
\end{equation*}

\noindent\textit{Proof of Theorem \ref{Thm:PurelySingularORAC}.}  Suppose $\mu_{ac}\neq 0$.  Then $\mu_{ac}$ must be $t\mu$, where $t$ is a nonzero scalar,   by Hutchinson's theorem.  Now, consider the unique Lebegue-Radon-Nikodym decomposition  (\ref{Eqn:Rudin6.10a})of $\mu$:
\begin{equation*}
\mu_{ac} = t\mu = t\mu_{ac} + t\mu_{s}.
\end{equation*}
But $\mu_{ac}\perp\mu_{s}$, so $(1 - t)\mu_{ac} = t\mu_{s}$ tells us that $\mu_{s} = 0$.  Therefore, $t = 1$ and $\mu_{ac} = \mu$.\hfill$\Box$

\begin{remark}\textbf{An alternate proof of Theorem \ref{Thm:PurelySingularORAC}.  }  Let $\widetilde{H}$ be a set operation defined by
\[ \widetilde{H}(S) = \bigcup_{b\in B} \tau_b(S).\]
By Hutchinson's theorem, the attractor $X$ of the IFS $(A,B)$ is the unique compact nonempty solution to $\widetilde{H}(X) = X$; the attractor $X$ is the support of the invariant measure $\mu$.  

Suppose $\mu_{ac}\neq 0$; then the support of $\mu_{ac}$ is not empty.  A quick calculation shows that
\begin{equation*}
\textrm{supp}(H\mu_{ac}) = \bigcup_{b\in B}\tau_b(\textrm{supp}(\mu_{ac})).
\end{equation*}   
However, we also know that $\mu_{ac} = H\mu_{ac}$ from (\ref{Eqn:MuACisHMuAC}).  Therefore,
\begin{equation*}
\textrm{supp}(\mu_{ac}) = \bigcup_{b\in B}\tau_b(\textrm{supp}(\mu_{ac})) = \widetilde{H}(\textrm{supp}(\mu_{ac})).
\end{equation*}
By the uniqueness for sets in Hutchinson's theorem, we can conclude that $\textrm{supp}(\mu_{ac}) = X$ and $\textrm{supp}(\mu_{s}) = \emptyset$.
\end{remark}
\begin{remark}\textbf{Relaxed conditions under which the dichotomy theorem is true.  } 

\noindent Suppose we work in an ambient space with measure $m$ and
\begin{enumerate}
\item $\mu$ is an \textit{equilibrium measure} for an IFS $S$ corresponding to a fixed
system of weights $\{p_i\}$---that is,
\begin{equation*}
\mu = \sum_{i} p_i \mu\circ\tau_i^{-1}
\end{equation*}
\item the maps $\tau_i$ scale $m$ by a positive constant---that is, there exist $a_i, b_i\in\mathbb{R}^{+}$ such that $a_i m \leq m\circ\tau_i^{-1}\leq b_i m$.
\end{enumerate}

If there is an analogue of the uniqueness theorem of Hutchinson for the equilibrium measure for $S$,  then the conclusions of Theorem \ref{Thm:PurelySingularORAC} still hold.
\end{remark}

\begin{example}Suppose $S$ is an IFS defined by positive weights and a rational map $z\mapsto r(z) = p(z)/q(z)$, where $p$ and $q$ are polynomials in one complex variable. The analogue of our $\tau$ maps (\ref{tauB}) are branches of the inverse of $r$, and the attractor $X$ is a Julia set \cite{Bea91}.  In this case, Brolin's theorem \cite{Bro65} guarantees the uniqueness of the equilibrium measure (up to scalar constant) $\mu$; $\mu$ is supported on $X$.  In this example, the measure $\mu$ satisfies a much more restrictive condition than $r$-invariance.  For all $\phi$, the measure $\mu$ satisfies
\begin{equation*}
\int \phi(z)d\mu(z) = \int \sum_{r(w) = z} p_w\phi(w) d\mu(z);
\end{equation*}
whereas $r$-invariance means that $\mu\circ r^{-1} = \mu$.
The equilibrium measure would be either absolutely continuous or purely singular and continuous with respect to the two-dimensional Lebesgue measure on the complex plane.  We note that these systems do not satisfy the conditions of Hutchinson's theorem.
\end{example}

\subsection{Gradations of the class of singular measures for affine IFSs}\label{Subsec:Gradation}

      The following argument for a particular class of measures $\mu$ on $\mathbb{R}^d$
shows that asymptotic properties of the Fourier transform are directly
related to how the measure $\mu$ translates. In particular the lower
bounds we obtain for the Fourier transform along particular paths to
infinity imply a discontinuity of the translates of $\mu$ near zero.

We state the next result for 1D, but note that it easily generalizes to $\mathbb{R}^d$ for $d > 1$.  The argument implies in particular the Riemann-Lebesgue conclusion for measures $\mu$ of bounded variation.  But it shows further that asymptotics of $\hat{\mu}$ are directly connected to the translation group $\mathbb{R}$ acting on the measures.  

For measures $\mu$ on $\mathbb{R}$, let $\|\mu\|$ denote the total variation norm.  For $t\in\mathbb{R}$ and Borel subsets $S\subset\mathbb{R}$, set
\begin{equation*}
(T_t(\mu))(S):=\mu(S-t)
\end{equation*}
Note that if $\mu = f\:dx$ for $f\in L^1(\mathbb{R})$ then $\|\mu\| = \|f\|_{L^1}$, and $\|f - T_t f\|_{L^1}\rightarrow 0$ as $t\rightarrow 0$.
 We get 
\begin{equation*}
\begin{split}
2\hat{\mu}(\xi)
& = \hat{\mu}(\xi) - e^{i\pi}\hat{\mu}(\xi)\\
& = \int e^{i 2\pi \xi x}\Biggl(d\mu(x) - d\mu\Bigl(x - \frac{1}{2\xi}\Bigr)\Biggr)\\
& = (\mu - T_{\frac{1}{2\xi}}\mu){\hat{\:\:}}(\xi);
\end{split}
\end{equation*}
and
\begin{equation}
|\hat{\mu}(\xi)| \leq \frac{1}{2} \|\mu -  T_{\frac{1}{2\xi}}\mu\|.
\end{equation}\label{Ineq:GradationProp}%
The next result follows:
\begin{proposition}
If $\{t_n\}\subset \mathbb{R} \backslash \{0\}$, then the implication $(i)\Rightarrow (ii)$ holds, where
\begin{enumerate}[(i)]
\item $\displaystyle \lim_{n\rightarrow\infty} \|\mu -  T_{\frac{1}{2t_n}}\mu\| = 0$
\item $\hat{\mu}(t_n)\rightarrow 0.$
\end{enumerate}
\end{proposition}
\textit{Proof:  }
Apply Inequality (\ref{Ineq:GradationProp}).
\hfill$\Box$
\begin{corollary}
Let $C\in \mathbb{R}_{+}$ be the constant in Theorem \ref{Thm:Erdos}, and let $\mu = \mu_{\lambda}$ be the measure from above, where $\lambda^{-1} = \alpha$ is a Pisot number.  Then
\begin{equation}
\| \mu - T_{\frac{\lambda^n}{2}}\mu \| \geq 2C
\end{equation}\label{Ineq:SeparationTranslates}%
for all $n\in\mathbb{N}$.
\end{corollary}\label{Cor:SeparationTranslates}

Fix $k\in\mathbb{N}$.  Following (\ref{Ineq:GradationProp}), we can use a factor of $e^{-i\pi/k}(1 - e^{i2\pi/k})$ instead of $2$ to obtain the following:
\begin{equation*}
|e^{-i\pi/k}(1 - e^{i2\pi/k})\hat{\mu}(\xi)|  = |(\mu - T_{\frac{1}{k\xi}}\mu){\hat{\:\:}}(\xi)|
\end{equation*}
or, more simply,
\begin{equation*}
\Big|2\sin\Bigl(\frac{\pi}{k}\Bigr) \hat{\mu}(\xi)\Big|= |(\mu - T_{\frac{1}{k\xi}}\mu){\hat{\:\:}}(\xi)|.
\end{equation*}
We have
\begin{equation*}
\Big|2\sin\Bigl(\frac{\pi}{k}\Bigr) \hat{\mu}(\xi)\Big|\leq \|(\mu - T_{\frac{1}{k\xi}}\mu){\hat{\:\:}}(\xi)\|.
\end{equation*}
For large $k$, $\sin(\frac{\pi}{k})\sim \frac{\pi}{k}$, so we can further grade the class of singular measures by asking which measures satisfy
\begin{equation*}
0 < \inf_{k,n} |k|\:\|\mu - T_{\frac{\lambda^n}{k}}\mu\|,
\end{equation*}
even if $\lambda$ is not the reciprocal of a Pisot number.  Finally, we could fix $\alpha$ and ask which measures satisfy
\begin{equation}
0 < \inf_{k,n} |k|^{\alpha}\:\|\mu - T_{\frac{\lambda^n}{k}}\mu\|.
\end{equation}\label{Ineq:AlphaGradation}

\noindent\textbf{Open problems:  }
\begin{enumerate}[(i)]
\item Give a geometric proof of a lower bound for the terms in (\ref{Ineq:SeparationTranslates}).
\item What is the closure of the set $\{T_{\frac{\lambda^n}{2}}(\mu_{\lambda})\:|\: n\in \mathbb{N}\}$?
\begin{remark}
Note that the set in (ii) is contained in the unit ball of the Banach space $BM$ of bounded-variation measures.  Since $BM$ is the dual of $C(X_{\lambda})$, the set in (ii) is relatively compact in the weak-star topology, so the closure in (ii) is compact, as is the closure of its convex hull $CV_{\lambda}$.
\end{remark}
\item What are the extreme points in $CV_{\lambda}$?
\item Find examples of measures which belong to one class in (\ref{Ineq:AlphaGradation}) but not another. 
\end{enumerate}

\subsection{Chaos and translates}\label{Subsec:Chaos}

Suppose $\{\tau_0, \tau_1\}$ is the IFS in (\ref{tauzeroone})  with invariant measure $\mu = \mu_{\lambda}$.  We ask what affine maps correspond to translation by $t$, where $t\in\mathbb{R}$.

\begin{lemma}  Let $B_0$ denote the set $\{0,1\}$.  If $\mu$ is the invariant measure of the affine IFS defined by $(\lambda, B_0)$, then $T_{t}\mu$ is the invariant measure of the affine IFS defined by $(\lambda,  B_t)$, where
\begin{equation}
B_t = \Bigl\{\frac{t}{\lambda},  \frac{t+\lambda}{\lambda}\Bigr\}.
\end{equation}\label{Eqn:Bt}%
\end{lemma}
\textit{Proof:  }
First, we consider the definition of $T_t\mu$: 
\begin{equation*}
T_t\mu(S) 
= \mu(S - t)
= \frac{1}{2}\mu\circ\tau_0^{-1}(S-t) + \frac{1}{2}\mu\circ\tau_1^{-1}(S-t).
\end{equation*}
We want to see which affine maps are associated with $\tau_0^{-1}\circ T_t$ and $\tau_1^{-1}\circ T_t$ where $T_t(s) = s - t$.
For each $s\in S$, 
\begin{equation*}
\tau_0^{-1}(s-t) = \lambda^{-1}(s-t) \quad \textrm{ and } \quad \tau_1^{-1}(s-t) = \lambda^{-1}(s-t) - 1;
\end{equation*}
writing $\tau_U^{-1}(s) = \lambda^{-1}(s-t)$, and $\tau_{V}^{-1}=\lambda^{-1}(s-t) - 1$ we find that the corresponding affine maps are 
\begin{equation*}
\tau_U(x) = \lambda\Bigl(x + \frac{t}{\lambda}\Bigr) \quad\textrm{ and }\quad \tau_V(x) = \lambda\Bigl(x + \frac{t+\lambda}{\lambda}\Bigr).
\end{equation*}
Therefore, the affine maps associated with translation have the same scaling factor $\lambda$  but different translations; the corresponding set $B_t$ is
\begin{equation*}
B_t = \Bigl\{\frac{t}{\lambda},  \frac{t+\lambda}{\lambda}\Bigr\}.
\end{equation*}
$ $\hfill$\Box$
Suppose $\lambda$ is the reciprocal of a Pisot number.  Now, in Corollary \ref{Cor:SeparationTranslates}, we considered $\|\mu - T_\frac{\lambda^n}{2}\mu\|$.  When we substitute $t = \frac{\lambda^n}{2}$ into $B_t$ in (\ref{Eqn:Bt}), we see that
\begin{equation*}
B_{\frac{\lambda^n}{2}} = \Biggl\{\frac{\frac{\lambda^n}{2}}{\lambda},  \frac{\frac{\lambda^n}{2}+\lambda}{\lambda}\Biggr\} = \Bigl\{ \frac{\lambda^{n-1}}{2},  \frac{\lambda^{n-1}}{2} + 1\Bigr\};
\end{equation*}
as $n\rightarrow\infty$, we see that the first element of $B_{\frac{\lambda^n}{2}}$ tends to $0$ and the second tends to $1$.  Loosely speaking, the two IFSs 
\[\{\tau_0, \tau_1\}\quad\textrm{ and }\quad\{\tau_{\frac{\lambda^{n-1}}{2}}, \tau_{\frac{\lambda^{n-1}}{2} + 1}\}\]
behave more and more like each other as $n$ increases, but their associated invariant measures stay apart by Corollary \ref{Cor:SeparationTranslates}.  We use this observation to motivate the following definition.

\begin{definition}
Suppose $B = \{b_i\}$ and let $B_t$ denote the set $\{b_i + At\} = \{b_i^{(t)}\}$.  Let $\mu^{(t)}$ denote the invariant measure of the system $(A, B_t)$.  We say that the measure $\mu^{(0)}:=\mu_{(A,B)}$ is \textit{chaotic} if and only if there exist a sequence $t_n\rightarrow 0$ and a positive number $\varepsilon$ such that
\begin{equation*}
\|\mu^{(t_n)} - \mu^{(0)}\| \geq \varepsilon
\end{equation*} 
for all $n$.
\end{definition}

\begin{theorem}
Suppose $\lambda$ is the inverse of a Pisot number and $(\lambda, B)$ is the system in (\ref{tauzeroone}).  Then the measure $\mu$ is chaotic.
\end{theorem}\label{Thm:ChaosErdos}%
\textit{Proof:  }Choose $\varepsilon = 2C$ in Corollary \ref{Cor:SeparationTranslates}, and choose 
\[t_n = \frac{\lambda^{n}}{2}.\]
We note that $k\in\mathbb{N}$, $k \geq 2$ could have been used instead of $2$ as in (\ref{Ineq:AlphaGradation}).
\hfill$\Box$

We call these measure ``chaotic'' because the phenomenon in Theorem \ref{Thm:ChaosErdos} mirrors that of chaos, in which two inputs of a system which are ``close'' to each other can result in drastically different outputs of that system.  In other words, the outcome of a small perturbation in initial conditions cannot be predicted.  Here, if the translation coefficients of an IFS are perturbed slightly, the total variation of the difference of the associated invariant measures will be separated by at least the constant in Corollary \ref{Cor:SeparationTranslates}.

Following this line of reasoning, we have three types of chaos in the IFSs in (\ref{generaltwotau}).
\begin{enumerate}
\item When $\lambda \in (0, 1/2)$, the systems are the ``most'' chaotic.  For example, when $\lambda = 1/3$, $\|\mu - T_{\frac{1}{k3^n}}\mu\| = 1$.
\item When $\lambda > 1/2$ is the reciprocal of a Pisot number, the systems are somewhat chaotic, by Corollary \ref{Cor:SeparationTranslates}.
\item For almost all $\lambda\in[1/2, 1)$, there is no chaos because the associated measures are absolutely continuous; the measures are represented by functions in $L^1(\mathbb{R})$, and translation is a continuous operation in $L^1$.
\end{enumerate}

Because the proof of \cite[Theorem 5.1]{LNR01} constructs a lower bound along a geometric sequence for the Fourier transform of the invariant measure for the affine IFS in $\mathbb{R}$, the measures in \cite[Theorem 5.1]{LNR01} are also chaotic.  See Remark \ref{Rem:RealExtensions}.

We note that when we translate IFSs in higher dimensions, the scaling matrix stays the same, but the translation coefficients shift, just as in $\mathbb{R}$.   If $(A,B)$ defines the IFS in (\ref{tauB}) with $b_i\in B$, the IFS corresponding to translating $\mu$ by $t\in\mathbb{R}^d$ will be defined by $(A,\widetilde{B})$, where $\tilde{b}_i = b_i + At \in \widetilde{B}$.  In this case, the measures in Theorems \ref{Thm:Erdos2D}, \ref{Thm:GenDir2D}, and \ref{Thm:ErdosdD} are also chaotic.

\section{Estimates for the Fourier transform in $\mathbb{R}^d$, $d = 2$}\label{Sec:2D}

In this section, we prove a two-dimensional analogue of Erd\H{o}s's theorem, where we restrict the Fourier transform of the Hutchinson measure to the line spanned by $[1,1]^t$ in $\mathbb{R}^2$.   We present the 2D proof before the general proof for $\mathbb{R}^d$ because it captures the essential ideas behind the proof of the higher-dimensional result, but the proof is simpler because it retains Erd\H{o}s's use of the cosine function after an initial estimate using the function $m_B$ (\ref{Eqn:GeneralmB}).  Then, we examine what happens when the entries of the direction vector are elements of $\mathbb{Z}^2$; in this case, we lose a bit of the flavor of Erd\H{o}s's proof but foreshadow techniques used in Section \ref{Sec:higherdim}.  We will finish the section by studying a measure associated with a matrix $A$ which is not a scalar multiple of the identity in Equation (\ref{Eqn:GeneralFTProduct}).

\subsection{A geometric progression in a fixed direction}\label{Subsec:FixedDirection2D}

Let $\lambda\in(0,1)$ be the inverse of a Pisot number $\alpha$, and let $\mathbf{x} \in \mathbb{R}^2$.  Using Equations (\ref{Eqn:GeneralmB}) and (\ref{Eqn:GeneralFTProduct}) with $A = \lambda^{-1}I$,
\begin{equation*}
B = \Biggl\{\begin{bmatrix}0\\0\end{bmatrix},\begin{bmatrix}1\\0\end{bmatrix},\begin{bmatrix}0\\1\end{bmatrix}\Biggr\},
\end{equation*}
and the equidistribution $p_b = \{1/3, 1/3, 1/3\}$, we eliminate the notation for dependence on $p_b$ and replace $A$ with $\lambda$ to write
\begin{equation*}
m_{B}(\mathbf{x}) = \frac{1}{3} \Bigl(1 + e^{i 2\pi x_1} + e^{i 2\pi x_2}\Bigr)
\end{equation*}
and
\begin{equation*}
\hat{\mu}_{\lambda,B}(\mathbf{x}) 
= 
\prod_{n = 1}^{\infty} m_{B}(\lambda^n \mathbf{x}).
\end{equation*}
When we restrict $\textbf{x}$ to the direction $W = \{\xi[1,1]^t : \xi\in\mathbb{R}\}$, the Fourier transform  of $\mu_{\lambda,B}$ becomes a function of a single variable $\xi$, and we denote the restricted Fourier transform of one variable as $\hat{\mu}_{\lambda,B,W}(\xi)$, where $\xi\in\mathbb{R}$:
\begin{equation}
\hat{\mu}_{\lambda, B, W}(\xi) 
=
 \prod_{n=1}^{\infty} \frac{1}{3}\Bigl(1 + 2\cos(2\pi \lambda^n \xi) 
+ i2\sin(2\pi \lambda^n \xi)\Bigr).
\end{equation}\label{Eqn:SingleVariable2D}%

\begin{theorem}
Let $\alpha$ be a Pisot number and let $\lambda = \alpha^{-1}$.  There exists a positive constant $C$ such that sequence $\{|\hat{\mu}_{\lambda, B, W}(\alpha^k)|\}_{k=0}^{\infty}$ is bounded from below by $C$; as a result, the sequence  $\{|\hat{\mu}_{\lambda, B, W}(\alpha^k)|\}_{k=0}^{\infty}$ does not tend to $0$ as $k\rightarrow\infty$.
\end{theorem}\label{Thm:Erdos2D}%
Before proving Theorem \ref{Thm:Erdos2D}, we state two immediate corollaries.
\begin{corollary}
If $\alpha$ is a Pisot number and $\lambda = \alpha^{-1}$, then the measure  $\mu_{\lambda}$ is not absolutely continuous with respect to Lesbegue measure on $\mathbb{R}^2$; in fact, $\mu_{\lambda}$ is purely continuous and singular.
\end{corollary}
\begin{corollary}
If $\alpha$ is a Pisot number and $\lambda = \alpha^{-1}$, then the measure  $\mu_{\lambda}$ is not the Fourier transform of a function $f\in L^1(\mathbb{R}^2)$.
\end{corollary}
\noindent\textit{Proof of Theorem \ref{Thm:Erdos2D}.  }Let $k\in\mathbb{N}$ and substitute $\xi = \alpha^k$ into Equation (\ref{Eqn:SingleVariable2D}) above.  We split $|\hat{\mu}_{\lambda, B, W}|^2$ into its real and imaginary parts to estimate
\begin{equation}
\begin{split}
|\hat{\mu}_{\lambda, B, W}(\alpha^k)|^2 
& = 
\prod_{n=1}^{\infty}
\Biggl(
\Bigl(\frac{1}{3} + \frac{2}{3}\cos(2\pi \lambda^n \alpha^k)\Bigr)^2 
+\Bigl(\frac{2}{3}\sin(2\pi \lambda^n \alpha^k)\Bigr)^2
\Biggr)\\
& \geq \prod_{n=1}^{\infty}\Bigl(\frac{1}{3} 
+ \frac{2}{3}\cos(2\pi \lambda^n \alpha^k)\Bigr)^2.
\end{split}
\end{equation}\label{Ineq:Cosine2D}%
Now, imitating Equations (\ref{Eqn:OneDProductPart1}) and (\ref{Eqn:OneDProductPart2}), the last product above can be written
\begin{equation}
|\hat{\mu}_{\lambda, B, W}(\alpha^k)|^2 \geq
\prod_{n=1}^{\infty}\Bigl(\frac{1}{3} + \frac{2}{3}\cos(2\pi \lambda^{n})\Bigr)^2
\prod_{n=0}^{k-1} \Bigl(\frac{1}{3} + \frac{2}{3}\cos(2\pi \alpha^{n})\Bigr)^2
\end{equation}\label{Eqn:ImitateErdos1D}%
 By Lemma \ref{Lemma:2DPositiveProduct} below, the infinite product in  (\ref{Eqn:ImitateErdos1D}) is positive.  

Choose $\theta\in (0,1)$ and $N\in\mathbb{N}$ such that two conditions are satisfied:
\begin{enumerate}[(a)]
\item the infinite product 
\begin{equation*}
\prod_{n=0}^{\infty} \Bigl(\frac{1}{3} + \frac{2}{3}\cos(2\pi \theta^{n})\Bigr)^2
\end{equation*}
is positive (Lemmas \ref{Lemma:Taylord} and \ref{Lemma:CountableTheta})
\item for all $n \geq N$, $\alpha^n < \theta^n < \frac{1}{4}$ mod $1$ (Lemma \ref{Lemma:GeometricBound}). 
\end{enumerate}

We now find ourselves back at the last part of Erd\H{o}s's theorem in Subsection \ref{Subsec:ErdosProof}. 

\bigskip

\noindent\textbf{Case 1:  }$k -1 \geq N$

\bigskip

 Define 
\begin{equation}
C = \prod_{n = 0}^{N-1} \Bigl(\frac{1}{3} + \frac{2}{3}\cos(2\pi \alpha^{n})\Bigr)^2;
\end{equation}\label{Eqn:Erdos2DC}%
we know $C>0$ by Lemma \ref{Lemma:FiniteProductNonzero2D}.  We have
\begin{equation*}
\begin{split}
|\hat{\mu}_{\lambda, B, W}(\alpha^k)|^2
& \geq \prod_{n=1}^{\infty}\Bigl(\frac{1}{3} + \frac{2}{3}\cos(2\pi \lambda^{n})\Bigr)^2
\prod_{n=0}^{k-1}\Bigl(\frac{1}{3} + \frac{2}{3}\cos(2\pi \alpha^{n})\Bigr)^2\\
& = C\prod_{n=1}^{\infty}\Bigl(\frac{1}{3} + \frac{2}{3}\cos(2\pi \lambda^{n})\Bigr)^2
\prod_{n=N}^{k-1} \Bigl(\frac{1}{3} + \frac{2}{3}\cos(2\pi \alpha^{n})\Bigr)^2\\
& \geq C\prod_{n=1}^{\infty}\Bigl(\frac{1}{3} + \frac{2}{3}\cos(2\pi \lambda^{n})\Bigr)^2
\prod_{n=N}^{k-1} \Bigl(\frac{1}{3} + \frac{2}{3}\cos(2\pi \theta^{n})\Bigr)^2\\
\end{split}
\end{equation*}
Just as in the proof of Erd\H{o}s's theorem in Subsection \ref{Subsec:ErdosProof}, we have found a lower bound which does not depend on $k$, since
\begin{equation*}
\prod_{n=N}^k \Bigl(\frac{1}{3} + \frac{2}{3}\cos(2\pi \theta^{n})\Bigr)^2 \geq \prod_{n=0}^\infty \Bigl(\frac{1}{3} + \frac{2}{3}\cos(2\pi \theta^{n})\Bigr)^2 > 0.
\end{equation*}%

\bigskip

\noindent\textbf{Case 2:  }$k \leq N$

\bigskip

If $k \leq N$, the finite product in (\ref{Eqn:ImitateErdos1D}) 
\begin{equation*}
\prod_{n=0}^k \Bigl(\frac{1}{3} + \frac{2}{3}\cos(2\pi \alpha^{n})\Bigr)^2
\end{equation*}
is not zero by Lemma \ref{Lemma:FiniteProductNonzero2D}.
\hfill $\Box$

\begin{lemma}  Let $\alpha$ be a Pisot number.  For any $j\in\mathbb{N}$, 
\begin{equation}
 \prod_{n = 0}^{j} \Bigl(\frac{1}{3} + \frac{2}{3}\cos(2\pi \alpha^{n})\Bigr)^2 > 0.
\end{equation}\label{Ineq:FiniteProductNonzero2D}%
\end{lemma}\label{Lemma:FiniteProductNonzero2D}%
\textit{Proof:  }
The product is nonzero because we know that
\begin{equation*}
\alpha^n \not\in \frac{1}{3} + \mathbb{Z} \text{ and } \alpha^n \not\in -\frac{1}{3} + \mathbb{Z}.
\end{equation*}
\hfill$\Box$
\begin{lemma}The product
\begin{equation*}
\prod_{n=1}^{\infty}\Bigl(\frac{1}{3} + \frac{2}{3}\cos(2\pi \lambda^{n})\Bigr)^2
\end{equation*}
is positive.
\end{lemma}\label{Lemma:2DPositiveProduct}%
\textit{Proof:  }
Using Lemma \ref{Lemma:Taylord} in the case $d = 2$, we can choose $N\in\mathbb{N}$  and $C > 0$ such that 
\begin{equation*}
\prod_{n=N}^{\infty} \Bigl(\frac{1}{3} + \frac{2}{3}\cos(2\pi \lambda^{n})\Bigr)^2 > C.
\end{equation*}
We now have to handle the first $N-1$ terms.  However, in order for for $1 \leq n \leq N-1$,
\begin{equation*}
\frac{1}{3} + \frac{2}{3}\cos(2\pi \lambda ^n) =  0,
\end{equation*}
we would need $\cos(2\pi \lambda^n) = -1/2$, which is only possible when
\begin{equation}
\lambda^n \in \frac{1}{3} + \mathbb{Z} \text{ or } \lambda^n \in -\frac{1}{3} + \mathbb{Z}.
\end{equation}\label{Eqn:RationalCondition}%
Since $\lambda$ is the inverse of a Pisot number, this is not possible.
\hfill$\Box$
We note that the proof of a higher-dimensional analogue of Theorem \ref{Thm:Erdos2D} (Theorem \ref{Thm:ErdosdD}) will be a bit more complicated because we are not guaranteed that $\arccos(-1/d)$ is a rational number multiplied by $2\pi$ as we have in the 2D case in Lemma \ref{Lemma:2DPositiveProduct}, (\ref{Eqn:RationalCondition}).  In order to get around this difficulty in Section \ref{Sec:higherdim}, we will make estimates using the function $m_B$ instead of with the cosine function.  We will see this same technique in Subsection \ref{Subsec:GeneralDirection2D}.

\subsection{General directions to infinity in the Fourier domain}\label{Subsec:GeneralDirection2D}

We retain the same notation from the previous section.

\begin{theorem}
Let $\alpha$ be a Pisot number and let $\lambda = \alpha^{-1}$.  Suppose $W = [n_1, n_2]^t\in \mathbb{Z}^2$.  There exists a positive lower bound for the sequence
$\Bigl\{|\hat{\mu}_{\lambda, B, W}(\alpha^k)|\Bigr\}_{k = 0}^{\infty}$.
\end{theorem}\label{Thm:GenDir2D}%

\textit{Proof:  }Without loss of generality, $n_1$ and $n_2$ can be positive integers. 

The following observation is important in what follows:  if $\textrm{dist}(\alpha^k, \mathbb{Z})\rightarrow 0$ as $k\rightarrow\infty$, and if $n\in\mathbb{Z}$, then
\begin{equation*}
\textrm{dist}(n\alpha^k, \mathbb{Z})\rightarrow 0 \text{ as } k\rightarrow\infty.
\end{equation*}
To see this, modify the proof of Lemma \ref{Lemma:Pisot}.
 
We have 
\begin{equation}
|\hat{\mu}_{\lambda, B, W}(\alpha^k)|^2 
 = \prod_{n=1}^{\infty}\frac{1}{9} \Big|1 + e^{i 2\pi n_1 \lambda^n \alpha^k} + e^{i2\pi  n_2 \lambda^n \alpha^k}\Big|^2.
\end{equation}\label{Eqn:SingleVariable2DIntegerDirection}%
We rewrite the product (\ref{Eqn:SingleVariable2DIntegerDirection}) as two products:
\begin{equation}
\prod_{n=1}^{\infty}\frac{1}{9} \Big|1 + e^{i 2\pi n_1 \lambda^n} + e^{ i2\pi  n_2 \lambda^n}\Big|^2
\end{equation}\label{Exp:InfiniteProductGeneral2D}%
and
\begin{equation}
\prod_{n=0}^{k-1}\frac{1}{9} \Big|1 + e^{i 2\pi n_1 \alpha^k} + e^{2\pi  n_2 \alpha^k}\Big|^2
\end{equation}
\label{Exp:FiniteProductGeneral2D}%

By a slight modification of Lemma \ref{Lemma:mBLlambda} in the case $d = 2$, we know that (\ref{Exp:InfiniteProductGeneral2D}) is nonzero.  See also (\ref{Exp:FiniteProduct2General2D}) below.

Now, we will focus on the product (\ref{Exp:FiniteProductGeneral2D}).  Choose $N\in\mathbb{N}$ and $\theta\in (0,1)$ such that for all $n \geq N$,
\begin{enumerate}[(i)]
\item $\displaystyle \prod_{n = 0}^{\infty}\frac{1}{9}\Bigl(1 + 2\cos(2\pi \theta^n)\Bigr)\neq 0$
\item  $n_1\alpha^n < \theta^n < \frac{1}{4} \text{ mod }1$
\item  $n_2\alpha^n < \theta^n < \frac{1}{4} \text{ mod }1$
\end{enumerate}

\bigskip

Suppose that $k-1 \geq N$.  To find a lower bound for the sequence 
$|\hat{\mu}_{\lambda, B, W}(\alpha^k)|^2$, we introduce the cosine terms into the estimate for the product (\ref{Exp:FiniteProductGeneral2D}):
\begin{equation*}
\begin{split}
&\prod_{n=0}^{k}\frac{1}{9} \Big|1 + e^{i 2\pi n_1 \alpha^n} + e^{ i2\pi  n_2 \alpha^n}\Big|^2\\
& \geq 
\prod_{n=0}^{N-1}\Big|1 + e^{i 2\pi n_1 \alpha^n} + e^{i 2\pi  n_2 \alpha^n}\Big|^2 \\
&\phantom{{\geq}} \cdot\prod_{n = N}^{k}\frac{1}{9}\Bigl(1 + \cos(2\pi n_1 \alpha^n) + \cos(2\pi n_2 \alpha^n)\Bigr)^2,
\end{split}
\end{equation*} 
For all $n \geq N$, all the terms $1$, $\cos(2\pi n_1 \alpha^n)$, and  $\cos(2\pi n_2 \alpha^n)$ are positive, and and we can say that
\begin{equation*}
\prod_{n = N}^{k}\frac{1}{9}\Bigl(1 + \cos(2\pi n_1 \alpha^n) + \cos(2\pi n_2 \alpha^n)\Bigr)^2
\geq  \prod_{n = N}^{k}\frac{1}{9}\Bigl(1 + 2\cos(2\pi \theta^n)\Bigr)^2.
\end{equation*}

In order to see if
\begin{equation}
\prod_{n=0}^{N-1}\Big|1 + e^{i 2\pi n_1 \alpha^n} + e^{i 2\pi  n_2 \alpha^n}\Big|^2 \\
\end{equation}\label{Exp:FiniteProduct2General2D}%
is nonzero, we consider how $1$, $e^{i2\pi t_1}$, and $e^{i2\pi t_2}$ can add to zero.  This is possible only when $t_1 = 1/3 (\text{ mod }1)$ and $t_2 = 2/3 (\text{ mod }1)$  or vice versa.  Therefore, since no integer multiple of $\alpha^n$ is a rational number, (\ref{Exp:FiniteProduct2General2D}) is nonzero.

\bigskip

If $k \leq N$, we do not need the cosine estimate.  By the same reasoning for (\ref{Exp:FiniteProduct2General2D}), we know that $|\hat{\mu}_{\lambda, B, W}(\alpha^k)|^2 \neq 0$.

\bigskip

If one of $n_1$ or $n_2$ is $0$, the discussion above simplifies.
\hfill$\Box$

\subsection{A family of Pisot matrices}
Suppose $\alpha$ is a Pisot number and suppose $c, b\in\mathbb{R}$ with $c > 1$.  Let $A$ be the $2\times 2$ matrix
\begin{equation}
A=\begin{bmatrix} \alpha & 0\\ b & c\end{bmatrix},
\end{equation}
and consider the measure defined by the affine IFS associated with $A$, $B = \{\textbf{0}, \textbf{e}_1, \textbf{e}_2\}$, and the equidistribution $p_b = \{1/3, 1/3, 1/3\}$.  By setting $\xi = \alpha^k[1, 0]^t$ in Equation (\ref{Eqn:GeneralFTProduct}), we obtain
\begin{equation*}
\hat{\mu}_{A,B,p}(\alpha^k [1,0]^t) = \prod_{n=1}^{\infty} m_B((A^t)^{-n}\alpha^k [1,0]^t).
\end{equation*} 
Now,
\begin{equation*}
(A^t)^{-n}\alpha^k [1,0]^t = [\alpha^{k-n}, 0]^t,
\end{equation*}
so
\begin{equation*}
\begin{split}
m_B([\alpha^{k-n}, 0]^t) 
& = \sum_{b\in B}\frac{1}{3}\Bigl(1 + e^{i 2\pi [1,0]\cdot [\alpha^{k-n}, 0]} + e^{i 2\pi [0,1]\cdot [\alpha^{k-n}, 0]}\Bigr)\\
&  = \frac{1}{3} (2 + e^{i 2\pi \alpha^{k-n}}).
\end{split}
\end{equation*}
At this point, we have reduced to the case $n_1 = 0$, $n_2 = 1$ in Subsection \ref{Subsec:GeneralDirection2D}, and we have a two-parameter family of $2\times 2$ matrices such that the sequence
\begin{equation*}
\{|\hat{\mu}_{A,B,p}(\alpha^k [1,0]^t)|\}
\end{equation*}
is bounded from below by a positive constant.

We note that the matrices here are not scalar multiples of orthonormal matrices, so the conclusion about singularity for the associated invariant measures does not follow from previous results in the literature, for example \cite[Theorem 3.1]{LNR01}.

\section{Estimates for the Fourier transform in $\mathbb{R}^d$, $d > 2$}\label{Sec:higherdim}

We return to the general case outlined in Subsection \ref{Subsec:HelpfulLemmas}:
\begin{equation*}
\begin{split}
 \hat{\mu}_{\lambda, B, W}(\xi) 
& = \prod_{n=1}^{\infty} \Biggl(\frac{1}{d+1}\Bigl(1 +d e^{2\pi i \lambda^n \xi}\Bigr)\Biggr)
= \prod_{n=1}^{\infty}  m_{B,W}(\lambda^n \xi)
\end{split}
\end{equation*}

Now, set $x = \alpha^k$, where $\alpha = \lambda^{-1}$.  We find a lower bound for $|\hat{\mu}_{\lambda, B, W}(\alpha^k)|^2$ for large $k$:  
\begin{equation*}
|\hat{\mu}_{\lambda, B, W}(\alpha^k)|^2
= \prod_{n=1}^{\infty} |m_{B,W}(\lambda^n \alpha^k)|^2
= \prod_{n=1}^{\infty}|m_{B,W}(\lambda^{n-k})|^2
\end{equation*}
Let $j = n-k$; when $n = 1$, $j = 1-k$.
\begin{equation*}
|\hat{\mu}_{\lambda, B, W}(\alpha^k)|^2
= \prod_{j+k=1}^{\infty}|m_{B,W}(\lambda^{j})|^2
=\prod_{j = 1-k}^{0}|m_{B,W}(\lambda^{j})|^2\prod_{j=1}^{\infty}|m_{B,W}(\lambda^{j})|^2
\end{equation*}
Finally, we can change the finite product into a function of $\alpha$:
\begin{equation}
|\hat{\mu}_{\lambda, B, W}(\alpha^k)|^2
=\prod_{n = 0}^{k-1}|m_{B,W}(\alpha^{n})|^2\prod_{n=1}^{\infty}|m_{B,W}(\lambda^{n})|^2
\end{equation}\label{Eqn:MuHatInitial}%

Set $C_{\lambda} = \prod_{n=1}^{\infty}|m_{B,W}(\lambda^{n})|^2$.  By Lemma \ref{Exp:FiniteProduct2General2D}, $C_{\lambda} > 0$.

With this notation,
\begin{equation*}
|\hat{\mu}_{\lambda, B, W}(\alpha^k)|^2
=  
C_{\lambda} \prod_{n = 0}^{k-1}|m_{B,W}(\alpha^{n})|^2.
\end{equation*}

By Lemmas \ref{Lemma:GeometricBound} and \ref{Lemma:CountableTheta}, we can now choose $\theta$ in $(0,1)$ and $N\in\mathbb{N}$ such that
\begin{enumerate}[(a)]
\item the following inequality is satisfied:
{\begin{equation}
 \prod_{n=0}^{\infty} \Bigl(\frac{1}{d+1}\Bigr)^2\Bigl(1 + d \cos(2\pi \theta^n) \Bigr)^2 >0
\end{equation}\label{necessaryfortheta}}
\item  the distance from $0$ to the
equivalence class of $\alpha^n$ in $\mathbb{R}/\mathbb{Z}$ is less
than $\theta^n < 1/4$ for each $n\geq N$.
\end{enumerate}

As in Subsection \ref{Subsec:ErdosProof}, we have chosen $N\in\mathbb{N}$ such that $\theta^n < 1/4$ for all $n \geq N$.  We make this choice so that for each $n \geq N$,
\begin{equation*} \cos(2\pi \alpha^n) \geq \cos(2\pi \theta^n).
\end{equation*}

\bigskip

\noindent\textbf{Case 1:  }$k - 1 \geq N$

\bigskip

For $k -1 \geq N$, we can write
\begin{equation*}
\begin{split}
|\hat{\mu}_{\lambda,B,W}(\alpha^k)|^2 = 
& C_{\lambda}
\prod_{n = 0}^{k-1}|m_{B,W}(\alpha^{n})|^2\\
&=
\prod_{n =0}^{N-1}|m_{B,W}(\alpha^{n})|^2\prod_{n =N}^{k-1}|m_{B,W}(\alpha^{n})|^2.
\end{split}
\end{equation*}
The constant $\prod_{n =0}^{N-1}|m_{B,W}(\alpha^{n})|^2$ is independent of $k$ and is nonzero because $d > 1$.  We now want to show that there exists $C > 0$ such that for all $k -1> N$,
\begin{equation*}
\prod_{n =N}^{k-1}|m_{B,W}(\alpha^{n})|^2 > C.
\end{equation*}
For each $n\in\mathbb{N}$,
\begin{equation*}
|m_{B,W}(\alpha^{n})|^2 
\geq \textrm{Re}\Bigl(|m_{B,W}(\alpha^{n})|\Bigr)^2 
= \Bigl(\frac{1}{d+1}\Bigr)^2\Bigl(1 + d\cos(2\pi \alpha^n)\Bigr)^2,
\end{equation*} 
and 
\begin{equation*} 
\prod_{n = N}^{k-1}|m_{B,W}(\alpha^{n})|^2
> \prod_{n = N}^{k-1} \Bigl(\frac{1}{d+1}\Bigr)^2\Bigl(1 + d\cos(2\pi \alpha^n)\Bigr)^2.
\end{equation*}
Finally, 
\begin{equation*}
\prod_{n = N}^{k-1} \Bigl(\frac{1}{d+1}\Bigr)^2\Bigl(1 + d\cos(2\pi \alpha^n)\Bigr)^2 
\geq \prod_{n = N}^{k-1} \Bigl(\frac{1}{d+1}\Bigr)^2\Bigl(1 + d\cos(2\pi \theta^n)\Bigr)^2,
\end{equation*}
and we can remove $k$ from the lower bound because
\begin{equation*}
\prod_{n = N}^{k-1} \Bigl(\frac{1}{d+1}\Bigr)^2\Bigl(1 + d\cos(2\pi \theta^n)\Bigr)^2 \geq \prod_{n = 0}^{\infty}\Bigl(\frac{1}{d+1}\Bigr)^2\Bigl(1 + d\cos(2\pi \theta^n)\Bigr)^2.
\end{equation*}

\bigskip

\noindent\textbf{Case 2:  }$k \leq N$

\bigskip

As before, we repeat our argument about the positivity of the finite product in \ref{Eqn:MuHatInitial}.

We have now proven 
\begin{theorem}
If $\lambda$ is the inverse of a Pisot number $\alpha$, then $\{|\hat{\mu}_{\lambda, B, W}(\alpha^k)|\}_{k=1}^{\infty}$ is bounded below by a positive constant. \end{theorem}\label{Thm:ErdosdD}%

\begin{corollary}The measure $\mu_{\lambda}$ is not absolutely continuous with respect to Lebesgue measure; $\mu_{\lambda}$ is purely continuous and singular.
\end{corollary}

\section{Induced measures}\label{Sec:InducedMeasures}

In this section, we make a systematic connection between
measures $\mu$ on infinite products $P$ and their induced measures $\nu$ on $\mathbb{R}^d$.  It is the Fourier transform $\hat{\nu}$ which we
study. Note that for a fixed measure on some compact infinite product
$P$, there are many induced measures on $\mathbb{R}^d$, in fact one
for each affine IFS.

In the models we consider for measures with compact support in
$\mathbb{R}^d$, there is an underlying \textit{coding space} $P$ and an \textit{encoding mapping} $\pi$. This will be made
precise in the present section. The essential ingredients in measures
from Hutchinson's theory \cite{Hut81} are a finite family of
contractive mappings $S$ with attractor $X(S)$ and an infinite
product measure $\mu$ defined on the Borel subsets of $P$. In a
wider context, we may consider infinite product spaces $P$ and
encoding mappings $\pi : P \rightarrow X$ where $X$ is compact in
$\mathbb{R}^d$, e.g., a Hutchinson attractor. For every measure $\mu$
on $P$ there is a pull-back (or induced) measure on $\mathbb{R}^d$ supported on $X$.  Generalizing Hutchinson's construction, we explore the analogous construction for determinantal measures on $P$
and their pull-backs to $\mathbb{R}^d$. Moreover we give a formula for the Fourier transform of the induced measures on $\mathbb{R}^d$.

The relevance of the determinantal measures is that they include important models from statistical mechanics and analysis, see, e.g., \cite{Jor06} and \cite{Lyo03}.   As we will see, determinantal measures are constructed from infinite matrices; when these matrices are non-diagonal, the off-diagonal entries capture correlations.  Determinantal measures are used in the analysis of infinite systems, especially in the study of long-range order.  Lemma \ref{Lemma:BorelFourier} shows the relevance of determinantal measures to asymptotics for Fourier transforms.

\subsection{Construction of determinantal measures}\label{Subsec:DetMeasureBackground}

If $T$ is an operator on $\ell_2$ satisfying Lemma \ref{Lemma:AppendixLemma}, we can associate with $T$ a determinantal measure $\mu = \mu_{T}$.  We use a matrix representation $(T_{i,j})$ of $T$ and matrix functions to define $\mu_T$ on cylinder sets in the infinite product space $P$ mentioned above.  The measures defined on the cylinder sets of $P$ will satisfy a Kolmogorov consistency condition (Lemma \ref{Lemma:ConsistencyRelations}), so the measure defined on cylinder sets will be well-defined on the Borel subsets of $P$ (Theorems \ref{Thm:Kolmogorov} and \ref{Thm:ExtendFrom1to2}).

Specifically, we start with the product space 
\begin{equation*}
P = \prod_{\mathbb{N}} \{0,1\} = \{0,1\}^{\mathbb{N}} = \{ \text{all functions }\omega:\mathbb{N}\rightarrow \{0,1\}\}
\end{equation*}
and an operator $T:\ell^2 \rightarrow \ell^2$ with matrix representation
\begin{equation*}
 T_{i,j}:=\langle \varepsilon_i, T\varepsilon_j\rangle
\end{equation*}
where $\{\varepsilon_{i}\}_{i\in\mathbb{N}}$ represents the canonical orthonormal basis (ONB) in $\ell^2$. Consider all finite subsets $F$ of $\mathbb{N}$ and the corresponding sets
\begin{equation*}
\begin{split}
\{0,1\}^F 
& = \prod_{F}\{0,1\} =  \{\text{all functions }\xi:F\rightarrow \{0,1\} \}\\
& = \{ \xi = (\xi_k) : k\in F,\: \xi_k\in \{0,1\} \}.\\
\end{split}
\end{equation*}
The finite sets $F$ allow us to use $(T_{i,j})$ to form $\#F\times\#F$
matrices as follows.  If $\xi\in \{0,1\}^F$, then the general cylinder set $G(\xi)$ is 
\begin{equation*}
\begin{split}
G(\xi) & :=\{ \omega\in P : \omega|_{F} = \xi\}\\ & = \{\omega\in P :
\omega_k = \xi_k \text{ for all }k\in F\}.
\end{split}
\end{equation*}
The measure $\mu_{T}$ is defined on the cylinder set $G(\xi)$ via
\begin{equation}
\mu_T(G(\xi)) = \det W(\xi),
\end{equation}\label{Eqn:MuAndW}%
where $W(\xi)$ is the following finite $\#F\times \#F$ matrix:
\begin{equation}
W(\xi)_{i,j}:=\Biggl( \xi_i \delta_{i,j} +
(-1)^{\xi_i}\Bigl(\delta_{i,j} - T_{i,j}\Bigr)\Biggr).
\end{equation}\label{Exp:WFmatrix}%
(See \cite[(7.7.9), p. 140]{Jor06}.)

To help see how the function $W$ works, we can consider $W$ on more
specialized cylinder sets.  This specialization will prove to be
useful by Lemma \ref{Lemma:OneSetsAreEnough}.  Again, let $F = (i_1, i_2,
\ldots, i_k)\subset \mathbb{N}$ be a finite set, and let $\xi^F$ have
the property $\xi_k^F = 1$ for all $k\in F$.  In this case,
$\mu_{T}(G(\xi^F))$ is just the determinant of the $\#F\times \#F$
submatrix of $(T_{i,j})$ which is formed by choosing rows $i_1,
\ldots, i_k$ and columns $i_1, \ldots, i_k$---that is,
\begin{equation}
\mu_{T}(G(\xi^F))
=\det\begin{pmatrix} 
T_{i_1, i_1} & T_{i_1, i_2} & \cdots & T_{i_1, i_k}\\ 
T_{i_2, i_1} & T_{i_2, i_2} & \cdots & T_{i_2, i_k}\\
\vdots       & \vdots       & \ddots & \vdots      \\
T_{i_k, i_1} & T_{i_k, i_2} & \cdots & T_{i_k, i_k}\\
\end{pmatrix}
.
\end{equation}\label{WF1submatrix}%
When $\xi_i = 1$ for all $i$ in
(\ref{Exp:WFmatrix}), then
$
\delta_{i,j} - (\delta_{i,j} - T_{i,j}) = T_{i,j}
$.  Even though the cylinder sets $G(\xi^F)$ are more specialized cylinder
sets, they actually determine the measure $\mu_T$, as explained in the Appendix.

\subsection{Determinantal measures and Fourier transforms}\label{Subsec:DetMeasure}

In this subsection, $P$ always refers to $\{0,1\}^{\mathbb{N}}$.

\begin{definition}
Let $p\in(0,1)$.  The $p$-Bernoulli measure on $P$ is the infinite product measure where each factor $1$ is assigned probability $p$ and each factor $0$ is assigned probability $q=1-p$.
\end{definition}\label{Def:pBernoulliMeasure}%
\begin{proposition}
Let $p\in(0,1)$.  Consider the determinantal measure $\mu_T$ for the $\infty\times\infty$ matrix
\begin{equation}
T_{i,j} = p(\delta_{i,j}),
\end{equation}\label{Exp:DiagonalT}%
i.e., $T$ is $p$ times the infinite identity matrix.  Then $\mu_T$ is the $p$-Bernoulli measure on $P$.
\end{proposition}\label{Prop:DiagonalT}%
\textit{Proof:  }
As already noted in Subsection \ref{Subsec:DetMeasureBackground},
$\mu_T$ is a probability measure on $P$ defined from a given
$\infty\times\infty$ matrix with spectrum in the interval $[0,1]$.

Let $F$ be a finite subset of $\mathbb{N}$ and let $\xi\in \{0,1\}^F$.  Substitution of (\ref{Exp:DiagonalT}) into (\ref{Exp:WFmatrix}) shows
that if $T$ is diagonal, then so is $W(\xi)$.  Specifically,
\begin{equation}
\begin{split}
W(\xi) 
&= \Biggl(\xi_i + (-1)^{\xi_i}(1-p)\Biggr)\delta_{i,j}\\
& = \begin{cases}
p\delta_{i,j} & \text{ if } \xi_i = 1\\
(1-p)\delta_{i,j} & \text{ if } \xi_i = 0.
\end{cases}
\end{split}
\end{equation}\label{Eqn:DiagonalW}%
Hence
\begin{equation*}
\mu_T(G(\xi)) = p^{\#\{\xi_i = 1\}}q^{\#\{\xi_i = 0\}},
\end{equation*}
which is the $p$-Bernoulli measure on $P$.
\hfill$\Box$
Later we will see how the $p$-Bernoulli measure on $P$ induces a familiar Hutchinson measure for the IFS (\ref{generaltwotau}) with $b_0= 0$ and $b_1 = \lambda^{-1}$.

The proof of the next result for Toeplitz matrices follows the same reasoning and will only be sketched.
\begin{proposition}
Let $a \in (0,1)$ be fixed, and set
\begin{equation}
T_{i,j} : = \frac{1-a}{1+a} a^{|i-j|}.
\end{equation}\label{Exp:TToeplitz}%
Then $T$ satisfies Lemma \ref{Lemma:AppendixLemma}.  If $F\subset \mathbb{N}$ is a finite subset, and if $\xi\in\{0,1\}^F$, then
\begin{equation}
W(\xi)_{i,i} = 
\begin{cases}
\frac{1-a}{1+a} & {\rm if\: } \xi_i=1\\
\frac{2a}{1+a}  & {\rm if\: } \xi_i=0
\end{cases};
\end{equation}\label{Eqn:WToeplitzDiagonal}%
while
\begin{equation}
W(\xi)_{i,j} = (2\xi_i - 1) \frac{1-a}{1+a} a^{|i-j|}
\end{equation}\label{Eqn:WToeplitzOffDiagonal}%
for $i\neq j$, i.e. the off-diagonal terms.
\end{proposition}\label{Prop:Toeplitz}%
\textit{Proof:  }
The result follows from a substitution of (\ref{Exp:TToeplitz}) into (\ref{Exp:WFmatrix}).  For the diagonal entries, we have
\begin{equation*}
W(\xi)_{i,i} = \xi_i + (-1)^{\xi_i}\Bigl(1 - \frac{1-a}{1+a}\Bigr),
\end{equation*}
and (\ref{Eqn:WToeplitzDiagonal}) follows.

For the off-diagonal entries, 
\begin{equation*}
W(\xi)_{i,j} = -(-1)^{\xi_i} \frac{1-a}{1+a} a^{|i-j|}.
\end{equation*}
Since $2\xi_i - 1 = -(-1)^{\xi_i}$, formula (\ref{Eqn:WToeplitzOffDiagonal}) follows.
\hfill$\Box$

We use the notation $\sigma_b,$ $b \in \{0,1\}$ to denote the right-shifts on $P$:
\begin{equation*}
\sigma_{b}(\omega_1\: \omega_2\: \cdots) :=(b\: \omega_1\: \omega_2\: \cdots).
\end{equation*}
Let $\tau_0$ and $\tau_1$  be contractive mappings from $\mathbb{R}^d$ into $\mathbb{R}^d$.  By \cite{Hut81}, there is a unique compact subset $X = X(\tau)\subset\mathbb{R}^d$ such that
\begin{equation}
X = \tau_0(X)\cup\tau_1(X).
\end{equation}\label{Exp:Attractor}%

\begin{lemma}
There is a unique continuous mapping 
\begin{equation*}
\pi:P \rightarrow X(\tau)
\end{equation*}
which is onto and which satisfies 
\begin{equation}
\pi \circ \sigma_b = \tau_b \circ\pi
\end{equation}\label{Eqn:SwitchPi}%
for $b\in\{0,1\}$.
\end{lemma}\label{Lemma:PiDefined}%
\textit{Proof:  }
For every $\omega \in P$, the following intersection
\begin{equation}
\bigcap_{n\in\mathbb{N}} \tau_{\omega|n} (X) = \{ \pi(\omega)\}
\end{equation}\label{Exp:Singleton}%
is a singleton, where
\begin{equation*}
\omega|n : = (\omega_1, \ldots, \omega_n),
\quad
\text{ and }
\quad
\tau_{\omega|n} := \tau_{\omega_1}\tau_{\omega_2} \cdots\tau_{\omega_n}.
\end{equation*}
The proof of the uniqueness assertion is left to the reader.
\hfill$\Box$

\begin{lemma}
For every Borel measure $\mu$ on $P$, set 
\begin{equation}
\nu : = \mu\circ\pi^{-1}
\quad
\text{ i.e., }
\quad  
\nu(E) = \mu(\pi^{-1}(E))
\end{equation}\label{Eqn:DefineNu}%
for Borel subsets $E\subset \mathbb{R}^d$, where
\begin{equation*}
\pi^{-1}(E) = \{ \omega \in P: \pi(\omega)\in E\}.
\end{equation*}
Then $\nu$ is a Borel measure on $\mathbb{R}^d$ supported on $X(\tau)$ with Fourier transform 
\begin{equation}
\hat{\nu}(t) = \int_P e^{i 2 \pi t\cdot \pi(\omega)} d\mu(\omega)
\end{equation}\label{Eqn:DetFourierTransform}%
for all $t = [t_1, t_2, \ldots, t_d]\in\mathbb{R}^d$.
\end{lemma}\label{Lemma:BorelFourier}%
\textit{Proof:  }
See \cite{JKS07a}.  Extending (\ref{Eqn:DefineNu}) we get the following transformation rule for integration
\begin{equation}
\int_{\mathbb{R}^d} f d(\mu\circ\pi^{-1}) = \int_P (f\circ\pi)d\mu
\end{equation}\label{Eqn:TransRule}%
for all $f$.  Indeed, the right-hand side in (\ref{Eqn:TransRule}) is a positive linear functional in $f$, and so by Riesz's theorem, the right-hand side defines integration with respect to a measure, which can be checked to be $\mu\circ\pi^{-1}$.

Setting $f(x) = e^{i2\pi t\cdot x}$ for the Fourier transform, we therefore have
\begin{equation*}
\hat{\nu}(t) = \int_{\mathbb{R}^d} e^{i 2\pi t\cdot x}d\mu\circ\pi^{-1}(x)
= \int_P e^{i2\pi t\cdot \pi(\omega)}d\mu(\omega).
\end{equation*}
$ $\hfill$\Box$
\begin{remark}The reader will notice that we used this formula in our derivation of our infinite-product representation for our Fourier transforms $\hat{\nu}(\cdot)$ in (\ref{Eqn:GeneralFTProduct}) and (\ref{Eqn:BernoulliConvolution}). 
\end{remark}
We recall that for the IFS (\ref{tauzeroone}) in $\mathbb{R}$, the shift $\sigma_0$ in $P$ corresponds to $\tau_0$ in $\mathbb{R}$, and the shift $\sigma_1$ corresponds to $\tau_1$ in $\mathbb{R}$.  With that in mind, we consider applying $\sigma_0$ and $\sigma_1$ to function $\xi:F\rightarrow\{0,1\}$.  If, for example, $\xi$ defines the cylinder set $(*, *, 1, *, 0, *, *, \ldots)$, that is $F = \{3, 5\}$ with $\xi(3) = 1$ and $\xi(5) = 0$, then $\sigma_i\xi$ is the cylinder set $(i, *, *, 1, *, 0, *, *, \ldots)$.  That is, we have a new function $\sigma_i\xi:\{1, 4, 6\}\rightarrow\{0,1\}$ such that $\sigma_i\xi(1) = i$, $\sigma_i\xi(4) = 1$, and $\sigma_i\xi(6) = 0$.
\begin{corollary}
Let $T=(T_{i,j})$ be an $\infty\times\infty$ matrix, $i,j \in \mathbb{N}$, with spectrum in $[0,1]$; let $W$ be the matrix function in (\ref{Exp:WFmatrix}); and let $\mu_{T}$ be the corresponding determinantal measure.

For $b\in\{0,1\}$ and $\omega = (\omega_1\:\omega_2\: \cdots)\in P$ set
\begin{equation*}
\sigma_{b}(\omega_1\: \omega_2\:\omega_3\: \cdots) :=(b\: \omega_1\: \omega_2\: \omega_3\: \cdots).
\end{equation*}
For all finite subsets $F\subset\mathbb{N}$ and all $\xi\in\{0,1\}^F$, we have the following recursive identity:
\begin{equation}
W(\sigma_0\xi) + W(\sigma_1\xi) 
= \begin{pmatrix}
1 & \:\:0\:\: & \:\:0\:\: & \:\:0\:\: & \cdots & 0\\
\hline \\
  &   &          &  &  &\\ 
2W'(\xi)  &   &  &  &&  \\
  &   &          &  &  &\\ 
\end{pmatrix}
\end{equation}\label{Eqn:RecursiveW}%
where $W'(\xi)$ is defined as in (\ref{Exp:WFmatrix}) but with respect to $T'_{i,j}:= T_{i+1,j}$.

In addition, suppose $\sum_{k}|T_{1,k}| < \infty$.  Then for the two measures $\mu_{T}\circ\sigma_b^{-1}$, $b\in\{0,1\}$ we have the following relative absolute continuity
\begin{equation}
\mu_{T}\circ \sigma_{b}^{-1} \ll \mu_{T}.
\end{equation}\label{Eqn:RelativeAbsoluteContinuity}%
\end{corollary}\label{Cor:ShiftMeasure}%
\textit{Proof:  }
The two conclusions follow from formulas (\ref{Eqn:MuAndW}) and (\ref{Exp:WFmatrix}) combined with basic determinant identities. 

From (\ref{Exp:WFmatrix}) it follows that the first rows in the matrices $W(\sigma_0\xi)$ and $W(\sigma_1\xi)$ from (\ref{Eqn:RecursiveW}) are
\begin{equation*}
( 1-T_{1,1}, -T_{1,2}, -T_{1,3}, \ldots)
\end{equation*}
and
\begin{equation*}
(T_{1,1}, T_{1,2}, T_{1,3}, \ldots).
\end{equation*} 
If $\sum_{k}|T_{1,k}| < \infty$, then the absolute continuity relations (\ref{Eqn:RelativeAbsoluteContinuity}) follow.
\hfill$\Box$

\subsection{Induced measures and Bernoulli IFSs in $\mathbb{R}$}

         We prove a general formula for the Fourier transform of the induced
measures for Bernoulli IFSs in 1D.

\begin{theorem}
Let $T:\ell^2\rightarrow \ell^2$ satisfy Lemma \ref{Lemma:AppendixLemma}, and let $\lambda\in (0,1)$ be given. 
Let $(T_{i,j})_{i,j\in\mathbb{N}}$ be the matrix representation; let $\mu_T$ be the determinantal measure on $P:=\{0,1\}^{\mathbb{N}}$.

Let $S_{\lambda}$ be the IFS from (\ref{generaltwotau}) with $b_0 = 0$ and $b_1 = \lambda^{-1}$.  Let 
\begin{equation}
\nu_{T} = \mu_T\circ\pi^{-1}
\end{equation}\label{Eqn:NuT}%
be the induced measure on $\mathbb{R}$.  For every $n\in\mathbb{N}$, set $F_n:=\{1,2,\ldots, n\}$, and let $T_{F_n}$ be the corresponding restricted matrix $(T_{i,j})_{i,j\in F_n}$.

Then the Fourier transform of $\nu_T$ satisfies
\begin{equation}
\hat{\nu}_T(t) = \lim_{n\rightarrow\infty} \det(I_n + D_n(\lambda t)T_{F_n})\text{ for all }t\in\mathbb{R}.
\end{equation}\label{Eqn:FourierNuT1}%
 The matrix $I_n$ is the $n\times n$ identity matrix, and $D_n(\lambda t)$ is defined by
\begin{equation}
D_n(\lambda t)  
= 
\begin{pmatrix}
e(\lambda t) - 1 & 0                 & \cdots & 0 \\
               0 & e(\lambda^2 t)-1  & \cdots & 0 \\
           \vdots& \vdots            & \ddots & \vdots \\
               0 & \cdots            &    0   & e(\lambda^n t)-1 \\
\end{pmatrix}
\end{equation}\label{Eqn:Dn}%
and
\begin{equation} e(x):=e^{i 2\pi x}, x\in\mathbb{R}.\end{equation}\label{Eqn:eofx}%
Moreover, we have the following asymptotic formula (referring to $n\rightarrow \infty$ in (\ref{Eqn:FourierNuT1}))
\begin{equation}
\hat{\nu}_T\:\substack{\simeq\\n\rightarrow\infty}
\:
 \exp\Biggl(\sum_{k=1}^n (e(\lambda^k t) - 1)T_{k,k}\Biggr)
\end{equation}\label{FourierNuT2}%
\end{theorem}\label{Thm:DetFourier1}%

\textit{Proof:  }
A computation shows that the encoding mapping $\pi:\{0,1\}^{\mathbb{N}}\rightarrow \mathbb{R}$ for the system $S_{\lambda}$ in (\ref{generaltwotau}) is
\begin{equation}
\pi_{\lambda}(\omega) = \sum_{k = 1}^{\infty}\omega_k\lambda^k
\end{equation}\label{Eqn:PowerSeriesPi}%
for $\omega = (\omega_1\:\omega_2\:\cdots)\in \{0,1\}^{\mathbb{N}}$.  Substituting (\ref{Eqn:PowerSeriesPi}) into (\ref{Eqn:DetFourierTransform}) in Lemma \ref{Lemma:BorelFourier}, and using Theorem \ref{Thm:ExtendFrom1to2}, we get the following limit formula:
\begin{equation}
\hat{\nu}_T(t)
= \lim_{n\rightarrow\infty}
\sum_{\omega\in \{0,1\}^{F_n}}
e\Bigl(\sum_{k = 1}^n \omega_k \lambda^k\Bigr)
\det W_{F_n}(\omega)
\end{equation}\label{Eqn:FourierNuT3}%
where $W_{F_n}$ refers to the $n\times n$ matrix
\begin{equation}
W_{F_n}(\omega) = \Bigl(\omega_i \delta_{i,j} + (-1)^{\omega_i}(\delta_{i,j} - T_{i,j}) \Bigr), \quad i,j\in F_n.
\end{equation}\label{Eqn:WFn}%

Theorem \ref{Thm:ExtendFrom1to2} and the estimates which follow justify an interchange of summation in (\ref{Eqn:FourierNuT3}).  Specifically, in carrying out the $\{0,1\}^{F_n}$-summation in (\ref{Eqn:FourierNuT3}), we may  do the individual sums
\begin{equation*}
\sum_{\omega_1\in\{0,1\}}, \sum_{\omega_2\in\{0,1\}}, \ldots 
\end{equation*}
one-by-one.  In the factorization
\begin{equation}
e\Bigl(\sum_{k = 1}^n \omega_k \lambda^k\Bigr) = \prod_{k = 1}^n e(\omega_k\lambda^k)
\end{equation}\label{Eqn:ExpSumToProduct}%
we may distribute the factors $e(\omega_k\lambda^k)$ on the rows in $W_{F_n}(\omega)$, $k = 1, 2, \ldots$.  For the $k$th row we get
\begin{equation}
v_k:=\delta_k + (e(\lambda^k t)-1)T_k
\end{equation}\label{Eqn:vk}%
where $\delta_k$ is a vector of all $0$s except for a $1$ in place $k$, and where
\begin{equation*}
T_k = (T_{k,j})_{j\in F_n}.
\end{equation*}
Hence the summation in (\ref{Eqn:FourierNuT3}) is the determinant of the $n\times n$ matrix
\begin{equation}
I_n + D_n(\lambda t)T_{F_n};
\end{equation}\label{Exp:IDTmatrix} %
the desired conclusion (\ref{Eqn:FourierNuT1}) follows.

An elementary result in matrix theory states that
\begin{equation}
\det(I+S)\simeq \exp(\textrm{trace}(S))
\end{equation}\label{Asymp:DetTrace}%
up to second order in $S$.  An application of this to $S = S_n = D_n(\lambda t)T_{F_n}$ for each $n$ yields the asymptotic formula (\ref{FourierNuT2}) in the statement of the theorem.  We can justify ignoring factors $S^2$ and higher since $\lambda < 1$ and
\begin{equation}
e(\lambda^k t) - 1 \sim \sin\Bigl(\frac{\lambda^k t}{2}\Bigr) \sim \frac{\lambda^k t}{2}
\end{equation}\label{Asymp:Sin}%
holds for $k$ sufficiently large, with error estimates governed by the terms in the Taylor expansion (Section \ref{Sec:higherdim}).
\hfill$\Box$

We note that in Theorem \ref{Thm:DetFourier1}, if the matrix $(T_{i,j})$ is diagonal, then we recover the product formula (\ref{Eqn:GeneralFTProduct}) with $A = \lambda$, $B = \{0,1/\lambda\}$, and probability distribution $\{1-p, p\}$.  In this case we get a slightly different set of maps from (\ref{tauzeroone}):
\begin{equation*}
\tau_0(x) = \lambda x  \qquad \tau_1(x) = \lambda x + 1.
\end{equation*}
The reader can check this directly using (\ref{Eqn:MuABp}) with the specified values of $A$ and $B$.  Also, when the matrix of $T$ is diagonal, the expression (\ref{Eqn:FourierNuT1}) involves only scalar functions.  However, in general, the analogue of  the infinite product (\ref{Eqn:GeneralFTProduct}), will involve more subtle matrix computations.   

\subsection{An example:  determinantal measures defined by Toeplitz matrices}
We will now give an explicit  formula for the Fourier transform (\ref{Eqn:DetFourierTransform}) when
specialized to the case when the given $\infty\times\infty$ matrix $(T_{i,j})$
is a Toeplitz matrix.  

Proposition \ref{Prop:DiagonalT} shows that if the matrix $(T_{i,j})$ for
$\mu_T$ is diagonal, we have the $p$-Bernoulli measure on $P$.  In turn, the $p$-Bernoulli measure on $P$ induces the Hutchinson measure on $X(S)$ (Section \ref{Sec:BackgroundIFS}), which we saw in Theorem \ref{Thm:DetFourier1}.  For
non-diagonal matrices, such as the Toeplitz matrices in Proposition
\ref{Prop:Toeplitz}, we get new and different IFS measures.

            We now outline the computation of the Fourier transform of
the induced measures, for the case when $T$ is Toeplitz. We show how
when $T$ is specialized to the determinantal measure (\ref{Exp:TToeplitz}) in Proposition \ref{Prop:Toeplitz} and the IFS is
specialized to our $\lambda$ system (\ref{generaltwotau}), we arrive at a
product formula for the Fourier transform of $\nu_T$, but our new formula differs
from the more familiar product formula
(\ref{Eqn:GeneralFTProduct}) for the Hutchinson
measures.

       Note that while the family of measures on $\mathbb{R}$ from Section
\ref{Sec:BackgroundIFS} above depends only on the single parameter $\lambda$,
for the new induced measures, there will be the additional parameter $a$
entering the definition (\ref{Exp:TToeplitz}) of the Toeplitz
matrix. So the induced measures in this case will have the pair
$(\lambda, a)$ as parameters.

\begin{corollary}
Let $\lambda\in (0,1)$ and let the IFS be as in (\ref{generaltwotau}).  Let the matrix $T$ in Theorem \ref{Thm:DetFourier1} be Toeplitz; i.e. pick real numbers $p,a\in (0,1)$ and set
\begin{equation}
T_{i,j} : = p a^{|i-j|}, \quad i,j\in\mathbb{N}.
\end{equation}\label{Eqn:ToeplitzAlphaA}%
(We could also use $\mathbb{Z}$ in the place of $\mathbb{N}$.)

Let $\mu_T = \mu_{p,a}$ be the corresponding determinantal measure on $P$, and let 
$\nu_{(\lambda,p,a)}
:=
\mu_{
(p,a)
}
\circ
\pi_{\lambda}^{-1}
$ be the induced measure on $\mathbb{R}$.  Then the $n$th approximation (see (\ref{Eqn:FourierNuT1})) to the Fourier transform $\hat{\nu}_{(\lambda,p,a)}(t)$ satisfies
\begin{equation}
\prod_{k=1}^n (p e(\lambda^kt) + 1 - p) + O(p^n).
\end{equation}\label{Exp:NthApproximationNuTHat}%
In particular, the Fourier asymptotics are the same for the diagonal matrix $T = (T_{i,j})= (p \delta_{i,j})$ and the $(p,a)$-Toeplitz matrix $T$ when applied to the Bernoulli system $S_{\lambda}$ in (\ref{tauzeroone}).   
\end{corollary}\label{Cor:ToeplitzFourier}%

Before beginning the proof, we note that  $\nu_{\lambda,p, a}$ is not absolutely continuous with respect to Lebesgue measure on $\mathbb{R}$ by Theorem \ref{Thm:Erdos}, but because $\nu_{\lambda,p, a}$ is not a Hutchinson measure defined by (\ref{Eqn:MuABp}), we cannot conclude that $\nu_{\lambda,p,a}$ is purely continuous and singular.
\textit{Proof:  }
In calculating the determinant of the $n\times n$ matrix in (\ref{Eqn:FourierNuT1}), use the fact that the anti-symmetric tensor space
\begin{equation*}
\wedge{}^n(\mathbb{C}^n) = \underbrace{\mathbb{C}^n\wedge\cdots\wedge\mathbb{C}^n}_{n\text{ times}}
\end{equation*}
is $1$-dimensional.  For $v_i\in\mathbb{C}^n$, $1 \leq i \leq n$, 
\begin{equation*}
v_1\wedge v_2\wedge \cdots \wedge v_n 
= \det
\begin{pmatrix}
v_1\\
v_2\\
\vdots\\
v_n\\
\end{pmatrix}
\boldsymbol{1}
\end{equation*}
where $\boldsymbol{1}$ denotes a unit-basis vector in $\wedge^n(\mathbb{C}^n)$.

Now substitute (\ref{Eqn:ToeplitzAlphaA}) into (\ref{Eqn:vk}).  For the row vector $T_k$ in (\ref{Eqn:vk}) we get
\begin{equation}
a^{k-1}\delta_1 + \cdots + a\delta_{k-1} + \delta_k + a\delta_{k+1} + \cdots + a^{n-k}\delta_n.
\end{equation}\label{Exp:RowKinT}%
Substituting into (\ref{Exp:RowKinT}) and using $\delta_k\wedge \delta_k = 0$, we get
\begin{equation*}
\begin{split}
& v_1\wedge v_2\wedge\cdots\wedge v_n\\
& = \prod_{k=1}^n (1 + p(e(\lambda^k t)-1))\:\delta_1\wedge\delta_2\wedge\cdots\wedge \delta_n\\
& \phantom{{==}} + p^n\prod_{k=1}^n(e(\lambda^kt)-1)\:T_1\wedge T_2\wedge\cdots\wedge T_n\\
& = \prod_{k=1}^n \Bigl(p e(\lambda^k t) + 1 - p \Bigr)\boldsymbol{1}
+ O\Bigl(p^n\prod_{k=1}^n(e(\lambda^kt)-1) \Bigr),
\end{split}
\end{equation*}
which is the desired conclusion.
\hfill$\Box$

We conclude with an exact formula for the Fourier transform of the measure
\begin{equation*}
\nu_{\lambda,p,a} := \mu_{T_{p,a}}\circ \pi_{\lambda}^{-1}.
\end{equation*}

\begin{lemma} 
Let 
\begin{equation*}A_n =
\begin{bmatrix}
1 & a & a^2 & \cdots &\cdots & a^{n-1}\\
a & 1 & a   & \cdots &\cdots & a^{n-2}\\
a^2 & a & 1 & \cdots&\cdots & a^{n-3}\\
\vdots & \vdots &\vdots & \ddots&  & \vdots\\
a^{n-2} & a^{n-3} &\cdots & a & 1 & a\\
a^{n-1} & a^{n-2} &\cdots  & a^2 & a & 1
\end{bmatrix}
\end{equation*}
Then $\det(A_n)
=
(1-a^2)^{n-1}$.
\end{lemma}\label{Lemma:An}%
\textit{Proof:  } 
By induction.
\hfill$\Box$

\begin{definition}
Denote by $T_n(\hat{k})$ the determinant of the $(n-1)\times(n-1)$ matrix obtained from $A_n$ by omitting the rows and columns at the $(k,k)$-place.  Further set $D(t) := e(t) - 1 = e^{i 2\pi t}-1$.
\end{definition}\label{Defn:TnHatK}%
For example, in Definition \ref{Defn:TnHatK}, when $n = 3$, $T(\hat{1}) = 1-a^2$, $T(\hat{2}) = 1-a^4$, and $T(\hat{3}) = 1-a^2$. 

\begin{corollary} We have
\begin{enumerate}[(a)]
\item
\begin{equation*}
\hat{\nu}_{\lambda,p,a}(t) = \lim_{n\rightarrow\infty} P_n(t)
\end{equation*}
where 
\begin{equation*}\begin{split}
P_n(t) 
& = 1 + p^{n-1}\sum_{k=1}^n T(\hat{k})\prod_{j\neq k} D(\lambda^j t) \\
& \phantom{{=1+}}+ p^n(1-a^2)^{n-1}\prod_{k=1}^n D(\lambda^k t),
\end{split}\end{equation*}
\item and when we pass to the limit $n\rightarrow\infty$, we get
\begin{equation*}
\hat{\nu}_{\lambda,p,a}(t) = 1 + \lim_{n\rightarrow\infty} \sum_{k=1}^n T_n(\hat{k}) \prod_{j\neq k} p D(\lambda^k t).
\end{equation*}
\end{enumerate}
\end{corollary}

\textit{Proof:  }
A simple Taylor series argument shows that the last term in (a) tends to $0$ as $n\rightarrow\infty$.  Our previous results justify the limit consideration for the second term in (a).

For the benefit of the reader we include some tensor considerations.  In the last formula in the proof of Corollary \ref{Cor:ToeplitzFourier}, we used the Grassmanian version $T_1\wedge T_2\wedge \cdots \wedge T_n$ for the determinant from Lemma \ref{Lemma:An}.  Recall the vectors $T_k$ (\ref{Exp:RowKinT}) are in $\mathbb{C}^n$, when $n$ is fixed.  Since convergence has already been established, the proof is complete.
\hfill$\Box$

\subsection{Infinite determinants}

Inspired by the determinants appearing in the previous subsections of Section \ref{Sec:InducedMeasures}, we define a function $\det_{\lambda}$ whose domain is the space of bounded self-adjoint operators on $\mathcal{H}$ and whose co-domain is the set of positive-definite functions on $\mathbb{R}$.  Recall the matrix $D_n(\lambda t)$ which is defined in (\ref{Eqn:Dn}).

\begin{definition}
Let $T$ be a bounded self-adjoint operator in an infinite-dimensional separable Hilbert space $\mathcal{H}$ which satisfies one of the equivalent conditions in Lemma \ref{Lemma:AppendixLemma}. Let $\lambda\in(0,1)$.  We define $\det_{\lambda}$ by
\begin{equation}
\rm det\it_{\lambda}(T) 
= \lim_{n\rightarrow\infty}
\det\Bigl((D_n(\lambda t)+I_n)T_{F_n} +( I_n -T_{F_n})\Bigr).
\end{equation}\label{Eqn:InfiniteDeterminant}%
\end{definition}\label{Defn:LambdaDet}%
\begin{corollary}
Let an operator $T$ be specified as in Definition \ref{Defn:LambdaDet}.
\begin{enumerate}[(i)]
\item The limit in (\ref{Eqn:InfiniteDeterminant}) exists pointwise.
\item The limit function $F_{\lambda, T}(t)$ is positive definite on $\mathbb{R}$ and continuous.
\item There is a measure $\nu_{\lambda,T}$ of compact support on $\mathbb{R}$ such that
\begin{equation}
F_{\lambda,T} = \hat{\nu}_{\lambda, T}(t)
\end{equation}\label{Eqn:FourierTransformFlambdaT}%
\item The function in (i) has a removable singularity at $t = 0$, with $F_{\lambda, T}(0)$ assigning the missing value for the limit expression at $t = 0$.
\item If $U:\mathcal{H}\rightarrow\mathcal{H}$ is a unitary operator leaving the (finite) linear combinations of $\{\varepsilon_i\}$ invariant, then
\begin{equation}
\det \hspace{0in}_{\lambda}(UTU^{*}) = \det\hspace{0in}_{\lambda}(T)
\end{equation}\label{Eqn:UnitaryDeterminant}%
for all $\lambda\in(0,1)$.
\end{enumerate}
\end{corollary}

\begin{remark}(Open problem)  Does the identity in (\ref{Eqn:UnitaryDeterminant}) hold for all unitary operators $U$?
\end{remark}

\appendix

\section{Determination of measures on infinite products}

It helps to identify three classes of subsets of infinite
products $P$. Specifically, here we take $P = \{0, 1\}^{\mathbb{N}}$,
or $P = \{0, 1\}^{\mathbb{Z}}$.  In the discussion below, we pick the
natural numbers $\mathbb{N} = \{1, 2, \ldots\}$ to be definite.

For the three classes, the first is contained in the second, and the
second in the third. The significance of the difference between the first
two classes and the third is that the first two may be prescribed by
certain \textbf{finite} configurations, while the extensions to class $3$ involves
transfinite induction.

The first two classes of subsets of $P$ are indexed by all finite
subsets of $\mathbb{N}$. Confusingly both classes are called cylinder
sets. Class $1$ consists of a special kind of cylinder set: if the
$1$s correspond to ``winnings'' and the $0$s correspond to ``losses,'' class $1$ amounts to specifying
a finite configuration of winnings.  In contrast, class $2$ specifies
finite configurations with prescribed winning and losing positions.
The idea is to get Borel measures on $P$ by extension from positive
set functions defined initially only on the special kinds of sets in
classes $1$ or $2$. Following Kolmogorov, we are concerned with two
assertions, existence and uniqueness.

\textit{Class $1$.  }   Pick $F$, an arbitrary finite subset in $\mathbb{N}$, and
consider a string of $1$s at the places identified by $F$, so a single
string of winnings is a point $w_F \in \{0,1\}^F$. Set $G(w_F)
= \{\omega \in P \::\: \omega|_F = w_F\}$.

\textit{Class $2$.  } Pick $F$ and consider all configurations (wins and losses)
$\xi \in \{0,1\}^F$. Set $G(\xi) = \{\omega \in P \::\: \omega|_F = \xi\}$.

\begin{definition}
The Borel $\sigma$-algebra (i.e., all the Borel subsets of
$P$) is the smallest $\sigma$-algebra of subsets of $P$ containing all the subsets
in class $2$.  Similarly the subsets in class $2$ form a basis for the compact
(Tychonoff) topology on $P$.
\end{definition}

\textit{Class $3$.  } All Borel subsets of $P$.

\begin{theorem}\rm(Kolmogorov).    
\it If a positive set function is defined on all
the sets in class 2 and is consistent on class 2, then it extends uniquely to all the
sets in class 3, hence to a positive Borel measure.
\end{theorem}\label{Thm:Kolmogorov} %

We now show that consistency holds for the sets in Class 2 with respect to the determinantal measures defined in Subsection \ref{Subsec:DetMeasureBackground}.
\begin{lemma}\rm (Consistency Relations)
\it
Let $T:\ell^2\rightarrow \ell^2$ be given and let $(T_{i,j})_{i,j\in\mathbb{N}}$ be the corresponding matrix representation.  Assume Lemma \ref{Lemma:AppendixLemma} holds.  Let $\mu_T$ be defined on the family of all cylinder subsets of $\{0,1\}^{\mathbb{N}}$ by formula (\ref{Eqn:MuAndW}). 

 Let $F\subset\mathbb{N}$ be a finite non-empty subset, and let 
$\xi\in\{0,1\}^F$.  For $k\in\mathbb{N}\backslash F$ define $\xi_{\pm}\in \{0,1\}^{F\cup \{k\}}$ as follows:
\begin{equation}
\xi_{+}(i)
= \begin{cases}
\xi(i) & \text{ if } i\in F\\
1      & \text{ if } i = k;
\end{cases}
\end{equation}
\begin{equation}
\xi_{-}(i)
= \begin{cases}
\xi(i) & \text{ if } i\in F\\
0      & \text{ if } i = k.
\end{cases}
\end{equation}

Then 
\begin{equation}
G(\xi_{+}) \cup G(\xi_{-}) = G(\xi),
\end{equation}\label{Eqn:FirstCon}%
\begin{equation}
G(\xi_{+}) \cap G(\xi_{-}) = \emptyset,
\end{equation}\label{Eqn:SecondCon}%
and
\begin{equation}
\mu_T (G(\xi_{+})) + \mu_T (G(\xi_{-})) = \mu_T(G(\xi)).
\end{equation}\label{Eqn:ThirdCon}%
\end{lemma}\label{Lemma:ConsistencyRelations}%
\textit{Proof:  }
Since formulas (\ref{Eqn:FirstCon}) and (\ref{Eqn:SecondCon}) are immediate from the definitions, we only need to check (\ref{Eqn:ThirdCon}).  Also note that by induction we can go from (\ref{Eqn:ThirdCon}) to stronger consistency relations, extending in steps from $\{0,1\}^F$ to $\{0,1\}^{\tilde{F}}$, where $\tilde{F}$ is any finite subset containing $F$.

In the extension from $F$ to $F\cup\{k\}$ we will be adding a row and a column to the matrix $W(\xi)$.  The position of the new row (and column) relative to the existing rows depends on where $k$ lies in the ordering of the points from $F$.  

In the row counting of the extended matrix, the new row is number $k$.  So, row $k$ in $W(\xi_{+})$ has entries $T_{k,j}$; while row $k$ in $W(\xi_{-})$ has entries $\delta_{k,j} - T_{k,j}$ with $j$ running over $F\cup\{k\}$.   Now consider the determinantal measure of the union $G(\xi_{+})\cup G(\xi_{-})$---this measure is $\mu_T(G(\xi_{+})) + \mu_T(G(\xi_{-}))$, which is the sum of determinants.  Using the fact that the determinant is multilinear, we can write $\mu_T(G(\xi_{+})) + \mu_T(G(\xi_{-}))$ as
\begin{equation}
\begin{matrix}   &          k                                      & \\
    \det         & \begin{bmatrix} W(\xi) & \vdots     & W (\xi)\\
                                        0\cdots 0 &  1     & 0\cdots 0 \\
                                   W(\xi) & \vdots     & W(\xi)\end{bmatrix} & k \\
\end{matrix}
\end{equation}\label{Eqn:DetRowColK}%
where we have indicated the position of row $k$ and column $k$ in the extended matrix.   Row $k$  in (\ref{Eqn:DetRowColK}) has the form $(0\: \cdots\: 0\:\: 1\:\: 0\:\cdots\: 0)$ with $1$ in the $k$th place.  But by algebra, the determinant in (\ref{Eqn:DetRowColK}) is
\begin{equation*}
\det(W(\xi)) = \mu_T(G(\xi))
\end{equation*}
and the proof of formula (\ref{Eqn:ThirdCon}) is completed.
\hfill$\Box$

Now, the construction in Subsection \ref{Subsec:DetMeasureBackground} will yield a Borel measure on $P$ if we know that we can extend the formula for determinants of submatrices (\ref{WF1submatrix}) uniquely to general cylinder sets in class 2 (\ref{Eqn:MuAndW}).  
\begin{theorem}\rm{(folklore, see e.g., \cite[Section 7.7]{Jor06}).  } \it
Let an $\infty\times\infty$ matrix $T$ be given as in Subsection \ref{Subsec:DetMeasureBackground}.

Define a positive set function $s$ on all the sets in class $1$ by
setting $s$ to be the determinant of the submatrix of $T$ determined
by $F$ as in (\ref{WF1submatrix}); then this function extends uniquely
to all the sets in class $2$ (and consistency holds).  By Kolmogorov,
therefore the function $s$ extends uniquely also to class $3$, hence to a
positive Borel measure.  This positive Borel measure is the
determinantal measure $\mu_{T}$.
\end{theorem}\label{Thm:ExtendFrom1to2}%
\noindent\textbf{Caution. } In general it is not true that a positive set function
defined on all the sets in class $1$, even with consistency on class $1$, extends to class
$2$, let alone to class $3$.

\begin{lemma}
The measures $\mu_{T}(G(\xi^F))$ of the cylinder sets defined with all
$1$s determine Equations (\ref{Eqn:MuAndW}) and (\ref{Exp:WFmatrix}),
thereby determining $\mu_{T}$ on all Borel subsets of $P$.
\end{lemma}\label{Lemma:OneSetsAreEnough}%
\textit{Proof:  }See Theorem \ref{Thm:ExtendFrom1to2}.\hfill$\Box$

\section{Facts from operator theory}
\begin{lemma}Let $T$ be a given bounded operator in a Hilbert space $\mathcal{H}$.  The following conditions are equivalent:
\begin{enumerate}[(i)]
\item$0 \leq \langle v|Tv\rangle \leq \|v\|^2$ for all $v\in\mathcal{H}$
\item There is an orthonormal basis (ONB) $\{\varepsilon_i\}$ for $\mathcal{H}$. such that the matrix $T_{i,j}:=\langle \varepsilon_i|T\varepsilon_j\rangle$ satisfies
\begin{equation}
0\leq \sum_i \sum_j \overline{c_i} \:T_{i,j}\:c_j \leq \sum_i |c_i|^2
\end{equation}\label{Eqn:A2}%
for all finite sequences $\{c_i\}$, $c_i\in\mathbb{C}$.
\item The property in (ii) holds for all ONBs.
\item For all finite subsets $F$ of the index set in (ii) we have
\begin{equation}
0 \leq \det(T_{i,j})_{i,j\in F} \leq 1
\end{equation}\label{Ineq:A3}%
where in (\ref{Ineq:A3}) we are taking the determinant of the submatrix from (\ref{Eqn:A2}) corresponding to rows and columns indexed by $F$.
\item $T = T^{*}$ and $\textrm{spectrum}(T)\subseteq [0,1]$.
\end{enumerate}
\end{lemma}\label{Lemma:AppendixLemma}%

\acknowledgements

The first named author has had helpful conversations with Dorin Dutkay, Sergei Silvestrov, and Oran Stenflo.  The third named author has had helpful conversations with Christopher French and David Romano.

\end{document}